\documentclass[10pt]{article}
\usepackage{bbm}
\usepackage{amssymb}
\usepackage{mathrsfs}
\usepackage{calc}
\usepackage{graphicx}
\usepackage[all,cmtip]{xy}
\usepackage{amsmath}
\usepackage{amscd}
\usepackage{lmodern}
\usepackage{amsfonts}
\usepackage{picinpar}
\usepackage{tikz-cd}
\usepackage{pgf,tikz,pgfplots}
\pgfplotsset{compat=1.14}
\usepackage{mathrsfs}
\usetikzlibrary{arrows}
\usepackage[amsmath,thmmarks]{ntheorem}
\usepackage{enumerate}
\usepackage{mathptmx}
\usepackage{stmaryrd}
\usepackage{array}
\usepackage{extarrows}
\usepackage[nottoc,numbib]{tocbibind}
\usepackage[]{hyperref}

\newtheorem{defi}{Definition}[section]
\newtheorem{thm}[defi]{Theorem}

\newtheorem{prop}[defi]{Proposition}
\newtheorem{cor}[defi]{Corollary}
\newtheorem{lem}[defi]{Lemma}
\theorembodyfont{}
\newtheorem{exa}[defi]{Example}
\theoremstyle{plain}
\newtheorem{rmk}[defi]{Remark}

\theoremstyle{nonumberplain} \theorembodyfont{} \theoremsymbol{\ensuremath{\square}} \qedsymbol{\ensuremath\square}
\newtheorem{pro}{Proof}

\DeclareMathOperator{\Gl}{GL}

\DeclareMathOperator{\Ker}{Ker}

\DeclareMathOperator{\Aut}{Aut}
\DeclareMathOperator{\Cent}{Cent}
\DeclareMathOperator{\Pic}{Pic}
\DeclareMathOperator{\Id}{Id}
\DeclareMathOperator{\Ele}{El}
\DeclareMathOperator{\Ind}{Ind}
\DeclareMathOperator{\Trace}{Tr}

\begin{document}
\definecolor{xdxdff}{rgb}{0.49019607843137253,0.49019607843137253,1}
\definecolor{ududff}{rgb}{0.30196078431372547,0.30196078431372547,1}

\title{Centralizers of elements of infinite order in plane Cremona groups}
\date{}
\author{ShengYuan Zhao}
\maketitle

\newcommand{\llpar}{(\!(}
\newcommand{\rrpar}{)\!)}
\newcommand{\GL}[2]{\operatorname{GL}_{#1}(\mathbf{#2})}
\newcommand{\PGL}[2]{\operatorname{PGL}_{#1}(\mathbf{#2})}
\newcommand{\PGLthrK}{\operatorname{PGL}_{3}(\mathbf{K})}
\newcommand{\AUK}{\operatorname{Aut}(\mathbb{P}_{\mathbf{K}}^2)}
\newcommand{\CRK}{\operatorname{Cr}_{2}(\mathbf{K})}
\newcommand{\CRC}{\operatorname{Cr}_{2}(\mathbf{C})}
\newcommand{\CRn}{\operatorname{Cr}_{n}(\mathbf{C})}
\newcommand{\charK}{\operatorname{char}(\mathbf{K})}
\newcommand{\Jonq}{\operatorname{Jonq}(\mathbf{K})}
\newcommand{\Jonqz}{\operatorname{Jonq}_0(\mathbf{K})}
\newcommand{\PGLtK}{\operatorname{PGL}_{2}(\mathbf{K})}
\newcommand{\CC}{\mathbf{C}}
\newcommand{\RR}{\mathbf{R}}
\newcommand{\NN}{\mathbf{N}}
\newcommand{\QQ}{\mathbf{Q}}
\newcommand{\ZZ}{\mathbf{Z}}
\newcommand{\FF}{\mathbb{F}}
\newcommand{\KK}{\mathbf{K}}
\newcommand{\Sym}{\mathscr{S}}
\newcommand{\HH}{\mathscr{H}}
\newcommand{\PP}{\mathbb{P}}
\newcommand{\AAA}{\mathbb{A}}
\newcommand{\SSS}{\mathbb{S}}
\newcommand{\Fpb}{\overline{\mathbf{F}_p}}
\newcommand{\Jonqui}{Jonqui\`eres }
\newcommand{\Kah}{K\"{a}hler}
\newcommand{\dashmapsto}{\mapstochar\dashrightarrow}

\newenvironment{ital}{\it}{}

\begin{abstract}
Let $\KK$ be an algebraically closed field. The Cremona group $\CRK$ is the group of birational transformations of the projective plane $\PP^2_{\KK}$.
We carry out an overall study of centralizers of elements of infinite order in $\CRK$ which leads to a classification of embeddings of $\ZZ^2$ into $\CRK$, as well as a classification of maximal non-torsion abelian subgroups of $\CRK$ when $\charK=0$. 

\end{abstract}
\renewcommand{\abstractname}{Abstract}

\setcounter{tocdepth}{2}
\tableofcontents

\section{Introduction}

Let $\KK$ be an algebraically closed field. The \emph{plane Cremona group} $\CRK$ is the group of birational transformations of the projective plane $\PP^2_{\KK}$. It is isomorphic to the group of $\KK$-algebra automorphisms of $\KK(X_1,X_2)$, the function field of $\PP^2_{\KK}$. Using a system of homogeneous coordinates $[x_0;x_1;x_2]$, a birational transformation $f\in\CRK$ can be written as
\[
[x_0:x_1:x_2]\dashmapsto [f_0(x_0,x_1,x_2):f_1(x_0,x_1,x_2):f_2(x_0,x_1,x_2)]
\]
where $f_0,f_1,f_2$ are homogeneous polynomials of the same degree without common factor. This degree does not depend on the system of homogeneous coordinates. We call it the \emph{degree} of $f$ and denote it by $deg(f)$. Geometrically it is the degree of the pull-back by $f$ of a general projective line. Birational transformations of degree $1$ are homographies and form $\AUK = \PGLthrK$, the group of automorphisms of the projective plane. 

\paragraph{Four types of elements.}
Following the work of M. Gizatullin, S. Cantat, J. Diller and C. Favre, we can classify an element $f\in\CRK$ into exactly one of the four following types according to the growth of the sequence $(deg(f^n))_{n\in \NN}$ (See \cite{DF01} Theorems 0.2, 0.3): 
\begin{enumerate}
	\item The sequence $(deg(f^n))_{n\in \NN}$ is bounded, $f$ is birationally conjugate to an automorphism $f'$ of a rational projective surface $X$ and a positive iterate of $f'$ lies in the connected component of the identity of the automorphism group $\Aut(X)$. We call $f$ an \emph{elliptic} element.
	\item The sequence $(deg(f^n))_{n\in \NN}$ grows linearly, $f$ preserves a unique pencil of rational curves and $f$ is not conjugate to an automorphism of any rational projective surface. We call $f$ a \emph{\Jonqui twist}.
	\item The sequence $(deg(f^n))_{n\in \NN}$ grows quadratically, $f$ is conjugate to an automorphism of a rational projective surface preserving a unique elliptic fibration. We call $f$ a \emph{Halphen twist}.
	\item The sequence $(deg(f^n))_{n\in \NN}$ grows exponentially and $f$ is called \emph{loxodromic}.
\end{enumerate}
If every element of a subgroup $G\subset \CRK$ is elliptic and if $G$ is conjugate to a group of automorphisms of a projective rational surface, then we call $G$ an \emph{elliptic group}.

The standard reference \cite{DF01} is written for $\KK=\CC$ but the same proof works over an algebraically closed field $\KK$ of characteristic different from $2$ and $3$. The only problem with characteristics $2$ and $3$ is that the important ingredient \cite{Giz80} does not deal with quasi-elliptic fibrations. This minor issue has been clarified in \cite{CanDol12} and \cite{CGL19} so that the above classification holds in arbitrary characteristic (cf. \cite{CanDol12} Section 4.1 and \cite{CGL19} Sections 4.3, 4.4). 

\paragraph{The \Jonqui group}
Fix an affine chart of $\PP^2$ with coordinates $(x,y)$.
\emph{The \Jonqui group} $\Jonq$ is the subgroup of the Cremona group of all transformations of the form
\begin{equation}\label{eq:jonq}
(x,y)\dashmapsto \left ( \frac{ax+b}{cx+d},\frac{A(x)y+B(x)}{C(x)y+D(x)} \right ) 
\end{equation}
where
\[
\begin{pmatrix}
	a&b\\c&d
\end{pmatrix}\in \PGLtK,\quad 
\begin{pmatrix}
	A&B\\C&D
\end{pmatrix}\in \operatorname{PGL}_2(\KK(x)).
\]
In other words, $\Jonq$ is the group of all birational transformations of $\PP^1\times\PP^1$ permuting the fibers of the projection onto the first factor; it is isomorphic to the semi-direct product $\PGLtK\ltimes \operatorname{PGL}_2(\KK(x))$. The \Jonqui group is defined only if a system of coordinates is chosen; a different choice of the affine chart yields a conjugation by an element of $\PGLthrK$. 
We can directly check with Formula \eqref{eq:jonq} that the sequence $(deg(f^n))_{n\in \NN}$ grows at most linearly. Thus elements of $\Jonq$ are either elliptic or \Jonqui twists. 
We denote by $\Jonqz$ the normal subgroup of $\Jonq$ that preserves fiberwise the rational fibration, i.e.\ the subgroup of those transformations of the form $(x,y)\dashmapsto \left ( x,\frac{A(x)y+B(x)}{C(x)y+D(x)} \right )$; it is isomorphic to $\operatorname{PGL}_2(\KK(x))$. A \Jonqui twist of the \Jonqui group will be called a \emph{base-wandering \Jonqui twist} if its action on the base of the rational fibration has infinite order. 

If $\KK=\Fpb$ is the algebraic closure of a finite field, then $\KK,\KK^*$ and $\PGLtK$ are all torsion groups. Thus, if $\KK=\Fpb$ then base-wandering \Jonqui twists do not exist. When $\charK=0$, or when $\charK=p>0$ and $\KK\neq\Fpb$, there exist base-wandering \Jonqui twists.

The group of automorphisms of a Hirzebruch surface will be systematically considered as a subgroup of the \Jonqui group in the following way:
\begin{equation}
\Aut(\FF_n)=\left\{(x,y)\dashmapsto \left(\frac{ax+b}{cx+d},\frac{y+t_0+t_1x+\cdots+t_nx^n}{(cx+d)^n}\right)\vert \begin{pmatrix}
	a&b\\c&d\end{pmatrix}\in \operatorname{GL}_{2}(\KK), t_0,\cdots,t_n\in\KK\right\}. \label{eq:autofhirz}
\end{equation}

\paragraph{Main results.}

\begin{thm}\label{zzthm}
Let $\Gamma$ be a subgroup of $\CRK$ which is isomorphic to $\ZZ^2$. Then there are a projective rational surface $X$ and a birational map $\phi:\PP^2\dashrightarrow X$ such that $\Gamma'=\phi\Gamma\phi^{-1}$ has a pair of generators $(f,g)$ that fits in one of the following (mutually exclusive) situations:
\begin{enumerate}
	\item $f,g$ are elliptic elements, $\Gamma'\subset \Aut(X)$.
	\item $f,g$ are Halphen twists preserving a same elliptic fibration on $X$, and $\Gamma'\subset \Aut(X)$.
	\item one or both of the $f,g$ are \Jonqui twists, and there exist $m,n\in\NN^*$ such that the finite index subgroup of $\Gamma$ generated by $f^m$ and $g^n$ is in an $1$-dimensional torus over $\KK(x)$ in $\Jonqz=\operatorname{PGL}_2(\KK(x))$; 
	\item $f$ is a base-wandering \Jonqui twist and $g$ is elliptic. In some affine chart, we can write $f,g$ in one of the following forms:
	\begin{itemize}
		\item $g$ is $(x,y)\mapsto(\alpha x,\beta y)$ and $f$ is $(x,y)\dashmapsto(\eta(x),yR(x^k))$ where $\alpha,\beta\in\KK^*, \alpha^k=1, R\in \KK(x),\eta\in \PGLtK, \eta(\alpha x)=\alpha \eta(x)$ and $\eta$ has infinite order;
		\item (only when $\charK=0$) $g$ is $(x,y)\mapsto(\alpha x,y+1)$ and $f$ is $(x,y)\dashmapsto(\eta(x),y+R(x))$ where $\alpha\in\KK^*, R\in \KK(x), R(\alpha x)=R(x), \eta\in \PGLtK, \eta(\alpha x)=\alpha \eta(x)$ and $\eta$ has infinite order.
	\end{itemize}
\end{enumerate}
\end{thm}
\begin{rmk}
When $\KK$ is the algebraic closure of a finite field, the above list can be shortened since there is no elliptic elements of infinite order nor base-wandering \Jonqui twists. 
\end{rmk}
\begin{rmk}\label{degreefunction}
From Theorem \ref{zzthm} it is not difficult to see that (see Section \ref{proofsection} for a proof), when $\Gamma$ is isomorphic to $\ZZ^2$, the degree function $deg:\Gamma\rightarrow \NN$ is governed by the word length function with respect to some generators in the following sense. In the first case of the above theorem it is bounded. In the second case it is up to a bounded term a positive definite quadratic form over $\ZZ^2$. In the third case the degree of $f^i\circ g^j$ is dominated by $deg(f) i+deg(g) j$ if $i,j$ are positive. In the fourth case the degree function is up to a bounded term $f^i\circ g^j\mapsto c \vert i\vert$ for some $c\in\QQ_+$. 
\end{rmk}


Theorem \ref{zzthm} is based on several known results. The main new feature is the following result (see Theorem \ref{mainthm} for a more precise reformulation). 
\begin{cor}\label{nozz}
Let $G\subset \Jonq$ be a subgroup isomorphic to $\ZZ^2$. Suppose that the action of $G$ on the base of the rational fibration is faithful. Then $G$ is an elliptic subgroup.
\end{cor}
We state it as a corollary (in cases 2, 3, 4 of Theorem \ref{zzthm} the action on the base is not faithful) but we will see that it is rather an intermediate step to prove Theorem \ref{zzthm}.

A \emph{maximal abelian subgroup} is an abelian subgroup which is not strictly contained in any other abelian subgroup. Over the field of complex numbers, finite abelian subgroups of $\CRC$ have been classified in \cite{Bla07} and maximal uncountable abelian subgroups of $\CRC$ have been studied in \cite{Des06a} (Theorems 1.5, 1.6). When $\charK=0$ we will use Theorem \ref{zzthm} to classify maximal abelian subgroups of $\CRK$ which contain at least one element of infinite order, see Theorem \ref{abmaxthm}. Our classification is more precise than the results in \cite{Des06a}. 

Another theorem we obtain from Theorem \ref{zzthm} is the following:
\begin{thm}\label{virtuallyabelian}
Assume $\charK=0$.
Let $f\in\CRK$ be an element of infinite order. If the centralizer of $f$ is not virtually abelian, then $f$ is an elliptic element and a power of $f$ is conjugate to an automorphism of $\AAA^2$ of the form $(x,y)\mapsto (x,y+1)$ or $(x,y)\mapsto (x,\beta y)$ with $\beta\in\KK^*$.
\end{thm}

\paragraph{Previously known results.}
Let us begin with the group of polynomial automorphism of the affine plane $\Aut(\AAA^2)$. It can be seen as a subgroup of $\CRK$. It is the amalgamated product of the group of affine automorphisms with the so called \emph{elementary group} 
\[\Ele(\KK)=\{(x,y)\mapsto (\alpha x +\beta, \gamma y+P(x))\vert \alpha,\beta,\gamma\in\KK, \alpha\beta\neq 0, P\in\KK[x]\}.\]
Let $\KK$ be the field of complex numbers. S. Friedland and J. Milnor showed in \cite{FM89} that an element of $\Aut(\AAA_{\CC}^2)$ is either conjugate to an element of $\Ele(\CC)$ or to a gengeralized H\'enon map, i.e.\ a composition $f_1\circ \cdots \circ f_n$ where the $f_i$ are H\'enon maps of the form $(x,y)\mapsto(y,P_i(y)-\delta_i x)$ with $\delta_i\in\CC^*$, $P_i\in\CC[y]$, $deg(P_i)\geq 2$. S. Lamy and C. Bisi showed in \cite{Lam01} and \cite{Bis04} that the centralizer in $\Aut(\CC^2)$ of a generalized H\'enon map is finite by cyclic, and that of an element of $\Ele(\CC)$ is uncountable (see also \cite{Bis08} for partial extensions to higher dimension). Note that, when viewed as elements of $\CRC$, a generalized H\'enon map is loxodromic and an element of $\Ele(\CC)$ is elliptic.

As regards the Cremona group, centralizers of loxodromic elements are known to be finite by cyclic (S. Cantat \cite{Can11} Theorem 5.1, J. Blanc-S. Cantat \cite{BC16} Corollary 4.7). Centralizers of Halphen twists are virtually abelian of rank at most $8$ (M. Gizatullin \cite{Giz80} Proposition 7, S. Cantat \cite{Can11} Proposition 5.2). When $\KK$ is the field of complex numbers, centralizers of elliptic elements of infinite order are completely described by J. Blanc-J. D\'eserti in \cite{BD15} Lemmas 2.7, 2.8 and centralizers of \Jonqui twists in $\Jonqz$ are completely described by D. Cerveau-J. D\'eserti in \cite{CD12}. Centralizers of base-wandering \Jonqui twists are also studied in \cite{CD12} but they were not fully understood, for example the results in loc.\ cit.\ are not sufficient for classifying pairs of \Jonqui twists generating a copy of $\ZZ^2$. Thus, in order to obtain a classification of embeddings of $\ZZ^2$ in $\CRK$, we need a detailed study of centralizers of base-wandering \Jonqui twists, which is the main task of this article. Regarding the elements of finite order and their centralizers in $\CRK$, the problem is of a rather different flavour and we refer the readers to \cite{Bla07}, \cite{DI09}, \cite{Ser10}, \cite{Ure20} and the references therein.

\begin{rmk}

There is a topology on $\CRK$, called Zariski toplogy, which is introduced by M. Demazure and J-P. Serre in \cite{Dem70} and \cite{Ser10}. Note that the Zariski topology does not make $\CRK$ an infinite dimensional algebraic group (cf. \cite{BF13}). With respect to the Zariski topology, the centralizer of any element of $\CRK$ is closed (J-P. Serre \cite{Ser10}). When $K$ is a local field, J. Blanc and J-P. Furter construct in \cite{BF13} an Euclidean topology on $\CRK$ which when restricted to $\PGLthrK$ coincides with the Euclidean topology of $\PGLthrK$; centralizers are also closed with respect to the Euclidean toplogy. In particular the intersection of the centralizer of an element in $\CRK$ with an algebraic subgroup $G$ of $\CRK$ is a closed subgroup of $G$, with respect to the Zariski topology of $G$ (and with respect to the Euclidean topology when the later is present).
\end{rmk}

\paragraph{Comparison with other results.}
S. Smale asked in the '60s if, in the group of diffeomorphisms of a compact manifold, the centralizer of a generic diffeomorphism consists only of its iterates. There has been a lot of work on this question, see for example \cite{BCW09} for an affirmative answer in the $C^1$ case. Similar phenomenons also appear in the group of germs of $1$-dimensional holomorphic diffeomorphisms at $0\in\CC$ (\cite{Eca81}). See the introduction of \cite{CD12} for more references in this direction. With regard to $\CRK$, it is known that loxodromic elements form a Zariski dense subset of $\CRK$ (cf. \cite{Xie15}, \cite{BedDil05}) and that their centralizers coincide with the cyclic group formed by their iterates up to finite index (cf. \cite{BC16}). Centralizers of general \Jonqui twists are also finite by cyclic (Remark \ref{exampledeserti}).
 
One may compare our classification of $\ZZ^2$ in $\CRK$ to the following two theorems where the situations are more rigid. The first can be seen as a continuous counterpart and is proved by F. Enriques \cite{Enriques} and M. Demazure \cite{Dem70}, the second can be seen as a torsion counterpart and is proved by A. Beauville \cite{Bea07}:
\begin{enumerate}
	\item \emph{If $\KK^{*r}$ embeds as an algebraic subgroup into $\CRK$, then $r\leq 2$; if $r=2$ then the embedding is conjugate to an embedding into the group of diagonal matrices $\Delta$ in $\operatorname{PGL}_3(\KK)$.}
	\item \emph{If $p\geq 5$ is a prime number different from the characteristic of $\KK$ and if $(\ZZ/p\ZZ)^{r}$ embeds into $\CRK$, then $r\leq 2$; if $r=2$ then the embedding is conjugate to an embedding into the group of diagonal matrices $\Delta$ in $\operatorname{PGL}_3(\KK)$.}
\end{enumerate}  
 
The classification of $\ZZ^2$ in $\CRK$ is a very natural special case of the study of finitely generated subgroups of $\CRK$; and information on centralizers can be useful for studying homomorphisms from other groups into $\CRK$, see for example \cite{Des06}. We refer the reader to the surveys \cite{Fav10},\cite{Can18} for representations of finitely generated groups into $\CRK$ and \cite{CX18} for general results in higher dimension. 

\paragraph{Acknowledgement.}
I would like to address my warmest thanks to my supervisor Serge Cantat for initiating me into Cremona groups, for numerous discussions, for his constant support and for encouraging me to write this paper. I would also like to thank JunYi Xie for helpful discussions on related topics ranging from proof details to general background. Many thanks to the anonymous referee for her/his careful reading and for her/his comments that improve largely the exposition of the paper.

\section{Elements which are not base-wandering \Jonqui twists}
This section contains a quick review of some scattered results about centralizers from \cite{Can11},\cite{BD15},\cite{CD12},\cite{BC16}. Some of the proofs are reproduced, because the original proofs were written over $\CC$ on the one hand, and because we will need some by-products of the proofs on the other hand.
\subsection{Loxodromic elements}
\begin{thm}[\cite{BC16} Corollary 4.7]\label{loxothm}
Let $f\in \CRK$ be a loxodromic element. The infinite cyclic group generated by $f$ is a finite index subgroup of the centralizer of $f$ in $\CRK$.
\end{thm}
\begin{pro}
We provide a proof which is simpler than \cite{BC16}. See \cite{Can11} and \cite{Can18} for the technical tools used in this proof. The Cremona group $\CRK$ acts faithfully by isometries on an infinite dimensional hyperbolic space $\mathbb{H}$ and the action of a loxodromic element is loxodromic in the sense of hyperbolic geometry. In particular there is a unique $f$-invariant geodesic $Ax(f)$ on which $f$ acts by translation and the translation length is $\log(\lim_{n\rightarrow \infty}deg(f^n)^{1/n})$. The centralizer $\Cent(f)$ preserves $Ax(f)$ and by considering translation lengths we get a morphism $\phi:\Cent(f)\rightarrow \RR$. We claim that the image of $\phi$ is discrete thus cyclic. Let us see first how the conclusion follows from the claim. Let $x\in \mathbb{H}$ be a point which corresponds to an ample class and let $y$ be an arbitrary point on $Ax(f)$. An element of the kernel $\Ker(\phi)$ fixes a point in $Ax(f)$ and thus fixes $Ax(f)$ pointwise because it commutes with $f$. Therefore for any element $g$ of $\Ker(\phi)$ the distance $d(x,g(x))$ is bounded by $2d(x,y)$.  This implies that $\Ker(\phi)$ is a subgroup of $\CRK$ of bounded degree. If $\Ker(\phi)$ were infinite then its Zariski closure $G$ in $\CRK$ would be an algebraic subgroup of strictly positive dimension contained, after conjugation, in the automorphism group of a rational surface. As $\Cent(f)$ is Zariski closed, the elements of $G$ commute with $f$. The orbits of a one-parameter subgroup of $G$ would form an $f$-invariant pencil of curves. This contradicts the fact that $f$ is loxodromic. Consequently $\Ker(\phi)$ is finite and hence $\Cent(f)$ is finite by cyclic.

Now let us prove the claim that the image of $\phi$ is discrete. This follows directly from a spectral gap property for translation lengths of loxodromic elements proved in \cite{BC16}. We give here an easier direct proof found with S. Cantat. Suppose by contradiction that there is a sequence $(g_n)_n$ of distinct elements of $\Cent(f)$ whose translation lengths on $Ax(f)$ tend to $0$ when $n$ goes to infinity. Without loss of generality, we can suppose the existence of a point $y$ on $Ax(f)$ and a real number $\epsilon>0$ such that $\forall n, d(y,g_n(y))<\epsilon$. Let $x\in\mathbb{H}$ be an element which corresponds to an ample class. Then it follows that \[\forall n, d(x,g_n(x))\leq d(x,y)+d(y,g_n(y))+d(g_n(y),g_n(x))<2d(x,y)+\epsilon=:d,\]
i.e.\ the sequence $(g_n)_n$ has bounded degree $d$. Up to extracting a subsequence, we can assume that all the $g_n$ has the same degree $k$ with $1<k\leq d$. Elements of degree $k$ of the Cremona group form a quasi-projective variety $\operatorname{Cr}_2^k(\KK)$. JunYi Xie proved in \cite{Xie15} that for any $0<\lambda<\log(k)$, the loxodromic elements of $\operatorname{Cr}_2^k(\KK)$ whose translation lengths are greater than $\lambda$ form a Zariski open dense subset of $\operatorname{Cr}_2^k(\KK)$. Thus the $g_n$ give rise to a strictly ascending chain of Zariski open subsets of $\operatorname{Cr}_2^k\color[rgb]{0,0,0}\color[rgb]{0,0,0}(\KK)$, contradicting the noetherian property of Zariski topology. This finishes the proof. Note that \cite{Xie15} is also used to prove the spectral gap property in \cite{BC16}.
\end{pro}

\subsection{Halphen twists}
We only recall here the final arguments of the proofs.
\begin{thm}[\cite{Giz80} and \cite{Can11} Proposition 4.7]\label{halphenthm}
Let $f\in \CRK$ be a Halphen twist. The centralizer $\Cent(f)$ of $f$ in $\CRK$ contains a finite index abelian subgroup of rank less than or equal to $8$. 
\end{thm}
\begin{pro}
Being a Halphen twist, the birational transformation $f$ is up to conjugation an automorphism of a projective rational surface and preserves a relatively minimal elliptic fibration. An iterate of $f$ preserves any fiber of this elliptic fibration (cf. \cite{Giz80} Proposition 7). Thus the $f$-invariant fibration is unique. As a consequence $\Cent(f)$ acts by automorphisms preserving this fibration. It is proved in \cite{Giz80} (see \cite{CGL19} Sections 4.3, 4.4 for a clarification in characteristics $2$ and $3$) that the automorphism group of a rational minimal elliptic surface has a finite index abelian subgroup of rank less than $8$.
\end{pro}

\subsection{Elliptic elements of infinite order}
In this section we reproduce a part of \cite{BD15}; we follow the original proofs (for $\charK=0$) in loc.\ cit.\ and some extra details are added in case $\charK>0$.

The proof of the following proposition in \cite{BD15} translates word by word to the case of postive characteristic. It is based on a $G$-Mori-program for rational surfaces due to J. Manin \cite{Man67} and V. Iskovskih \cite{Isk79}.
\begin{prop}[\cite{BD15} Proposition 2.1]\label{GMori}
Let $S$ be a smooth rational surface over $\KK$. Let $f\in \Aut (S)$ be an automorphism of infinite order whose action on $\Pic(S)$ is of finite order. Then there exists a birational morphism $S\rightarrow X$ where $X$ is a Hirzebruch surface $\FF_n$ ($n\neq 1$) or the projective plane $\PP^2$, which conjugates $f$ to an automorphism of $X$.
\end{prop}

\begin{prop}[\cite{BD15} Proposition 2.3]\label{ellipticnormalform}
Let $f\in \CRK$ be an elliptic element of infinite order. Then $f$ is conjugate to an automorphism of $\PP^2$. Furthermore there exists an affine chart with affine coordinates $(x,y)$ on which $f$ acts by automorphism of the following form:
\begin{enumerate}
	\item $(x,y)\mapsto (\alpha x,\beta y)$ where $\alpha,\beta \in \KK^*$ are such that the kernel of the group homomorphism $\ZZ^2\rightarrow \KK^*,(i,j)\mapsto\alpha^i\beta^j$ is generated by $(k,0)$ for some $k\in\ZZ$;
	\item $(x,y)\mapsto (\alpha x,y+1)$ where $\alpha\in \KK^*$ and $\alpha$ is of infinite order if $\charK>0$.
\end{enumerate}
\end{prop}
\begin{rmk}
If $\KK=\Fpb$ then every elliptic element is of finite order.
\end{rmk}

\begin{pro}[of Proposition \ref{ellipticnormalform}]
By Proposition \ref{GMori} we can suppose that $f$ is an automorphism of $\PP^2$ or of a Hirzebruch surface.

Let's consider first the case when $f\in\Aut(\PP^2)=\PGLthrK$. By putting the corresponding matrix in Jordan normal form, we can find an affine chart on which $f$ has, up to conjugation, one of the following forms:
1) $(x,y)\mapsto (\alpha x,\beta y)$; 2) $(x,y)\mapsto (\alpha x,y+1)$; 3) $(x,y)\mapsto (x+y,y+1)$.
If $\charK >0$ then $f$ can not have the third form since it would have finite order; if $\charK =0$ then in the third case $f$ is conjugate by $(x,y)\dashmapsto (x-\frac{1}{2}y(y-1),y)$ to $(x,y)\mapsto (x,y+1)$. We now show that in the first case $\alpha,\beta$ can be chosen to verify the conditon in the proposition. Let $\phi:(x,y)\mapsto (\alpha x,\beta y)$ be a diagonal automorphism, we denote by $\Delta(\phi)$ the kernel of the group morphism $\ZZ^2\rightarrow \KK^*,(i,j)\mapsto \alpha^i\beta^j$. For $M=\begin{pmatrix}a&b\\c&d\end{pmatrix}\in \GL{2}{Z}$, we denote by $M(\phi)$ the diagonal automorphism $(x,y)\mapsto (\alpha^a\beta^bx,\alpha^c\beta^dy)$, i.e.\ the conjugate of $\phi$ by the monomial map $(x,y)\dashmapsto(x^ay^b,x^cy^d)$.
We have the relation $\Delta(M(\phi))=(M^{\intercal})^{-1}(\Delta(\phi))$. Using the Smith normal form of a matrix with integer entries, up to conjugation by a monomial map we can suppose that our elliptic element $f$ satisfies $\Delta(f)=<(k_1,0),(0,k_1k_2)>$ where $k_1,k_2\in \ZZ$. Since $f$ is of infinite order, $k_1k_2$ must be $0$.

If $f\in\Aut(\FF_0)=\Aut(\PP^1\times\PP^1)$, then we reduce to the case of $\PP^2$ by blowing up a fixed point and contracting the strict transforms of the two rulings passing through the point. If $f\in\Aut(\FF_n)$ for $n\geq 2$ and if $f$ has a fixed point which is not on the exceptional section, then we can reduce to $\FF_{n-1}$ by making an elementary transformation at the fixed point.

Suppose now that $f\in\Aut(\FF_n), n\geq 2$ and its fixed points are all on the exceptional section. By removing the exceptional section and an invariant fiber of the rational fibration, we get an open subset isomorphic to $\AAA ^2$ on which $f$ can be written as:
$(x,y)\mapsto(\alpha x,\beta y+Q(x))$ or $(x,y)\mapsto(x+1,\beta y+Q(x))$ where $\alpha,\beta\in \KK^*$ and $Q$ is a polynomial of degree $\leq n$. 

In the first case, the fact that there is no extra fixed point on the fiber $x=0$ implies $\beta=1$ and $Q(0)\neq 0$. The action on the fiber at infinity can be obtained by a change of variables $(x',y')=(1/x,y/x^n)$, so the fact that there is no extra fixed point on it implies $\beta=\alpha^n$ and $\deg(Q)=n$. This forces $\alpha$ to be a primitive $r$-th root of unity for some $r\in\NN$. Conjugating $f$ by $(x,y)\mapsto(x,y+\gamma x^d)$, we replace $Q(x)$ with $Q(x)+\gamma(\alpha^d-1)x^d$. This allows us to eliminate the term $x^d$ of $Q$ unless $\alpha^d=1$. So we can assume that $f$ has the form $(x,y)\mapsto (\alpha x, y+\tilde{Q}(x^r))$ where $\alpha^r=1$ and $\tilde{Q}\in\KK[x]$. Then $f$ is conjugate to $(x,y)\mapsto (\alpha x,y+1)$ by $(x,y)\dashmapsto(x,y/\tilde{Q}(x^r))$. Remark that this case does not happen in positive characteristic because an automorphism of this form would have finite order. Note that in this paragraph we did not use the fact that $f$ is of infinite order, so that Proposition \ref{jonqzelliptic} is proved.

Suppose now we are in the second case. There is no extra fixed point if and only if $\beta=1$ and $\deg (Q)=n$. If $\charK>0$ and if $\beta=1$, then $f$ would be of finite order. Therefore we can assume $\charK=0$. In that case, we can decrease the degree of $Q$ by conjugating $f$ by a well chosen birational transformation of the form $(x,y)\dashmapsto(x,y+\gamma x^{n+1})$ with $\gamma\in \KK^*$. By induction we get $(x,y)\mapsto (x+1,y)$ at last.
\end{pro}

Applying the proof of Proposition \ref{ellipticnormalform} to an automorphism preserving the fibration fiber by fiber, we get the following:
\begin{prop}\label{jonqzelliptic}
Let $f$ be an automorphism of a Hirzebruch surface which preserves the rational fibration fiber by fiber (we do not assume that $f$ is of infinite order). Then there exists an affine chart on which $f$ acts as an automorphism of the following form:
\begin{enumerate}
	\item $(x,y)\mapsto (x,\beta y)$ where $\beta \in \KK^*$;
	\item $(x,y)\mapsto (x,y+1)$.
\end{enumerate}
Here $x$ is the coordinate on the base of the rational fibration.
\end{prop}

Once we have the above normal forms, explicit calculations can be done:
\begin{thm}[\cite{BD15} Lemmas 2.7 and 2.8]\label{ellipticthm}
Let $f\in\CRK$ be an elliptic element of infinite order. 
\begin{enumerate}
	\item If $f$ is of the form $(x,y)\mapsto(\alpha x,\beta y)$ where $\alpha,\beta \in \KK^*$ are such that the kernel of the group homomorphism $\ZZ^2\rightarrow \KK^*,(i,j)\mapsto\alpha^i\beta^j$ is generated by $(k,0)$ for some $k\in\ZZ$, then the centralizer of $f$ in $\CRK$ is
	\[
	\Cent(f)=\{(x,y)\dashmapsto(\eta(x),yR(x^k))\vert R\in \KK(x),\eta\in \PGLtK, \eta(\alpha x)=\alpha \eta(x)\}.
	\]
	\item If $\charK=0$ and if $f$ is of the form $(x,y)\mapsto(\alpha x,y+1)$, then the centralizer of $f$ in $\CRK$ is
	\[
	\Cent(f)=\{(x,y)\dashmapsto(\eta(x),y+R(x))\vert \eta\in \PGLtK, \eta(\alpha x)=\alpha \eta(x), R\in \KK(x), R(\alpha x)=R(x)\}.
	\]
	If $\alpha$ is not a root of unity then $R$ must be constant and $\eta(x)=\beta x$ for some $\beta\in\KK^*$.
	\item If $\charK=p>0$ and if $f$ is of the form $(x,y)\mapsto(\alpha x,y+1)$ (where $\alpha$ must be of infinite order), then the centralizer of $f$ in $\CRK$ is
	\[
	\Cent(f)=\{(x,y)\dashmapsto(R(y)x,y+t)\vert t\in\KK, R(y)=S(y)S(y-1)\cdots S(y-p+1), S\in \KK(y)\}.
	\]
\end{enumerate}
\end{thm}
\begin{rmk}\label{pglremark}
Note that the elements $\eta$ in the above statement can be explicitly described:
\[
\{\eta\in\PGLtK\vert\eta(\alpha x)=\alpha \eta(x)\}=\begin{cases}\PGLtK\quad \text{if} \quad \alpha=1\\
\{x\mapsto\gamma x^{\pm 1}\vert \gamma\in\KK^*\}\quad \text{if} \quad \alpha=-1\\
\{x\mapsto\gamma x\vert \gamma\in\KK^*\}\quad \text{if} \quad \alpha\neq\pm 1\\
\end{cases}
\]
\end{rmk}

\begin{pro}
\emph{First case.} We treat first the case where $f$ is of the form $(x,y)\mapsto(\alpha x,\beta y)$. Let $(x,y)\dashmapsto(\frac{P_1(x,y)}{Q_1(x,y)},\frac{P_2(x,y)}{Q_2(x,y)})$ be an element of $\Cent(f)$; here $P_1,P_2,Q_1,Q_2\in\KK[x,y]$. The commutation relation gives us 
\[\frac{P_1(\alpha x,\beta y)}{Q_1(\alpha x,\beta y)}=\frac{\alpha P_1(x,y)}{Q_1(x,y)},\quad \frac{P_2(\alpha x,\beta y)}{Q_2(\alpha x,\beta y)}=\frac{\beta P_2(x,y)}{Q_2(x,y)}\]
which imply that $P_1,P_2,Q_1,Q_2$ are eigenvectors of the $\KK$-linear automorphism $\KK[x,y]\rightarrow\KK[x,y], g(x,y)\mapsto g(\alpha x,\beta y)$. Therefore each one of the $P_1,P_2,Q_1,Q_2$ is a product of a monomial in $x,y$ with a polynomial in $\KK[x^k]$.
Then we must have $\frac{P_1(x,y)}{Q_1(x,y)}=xR_1(x^k)$ and $\frac{P_2(x,y)}{Q_2(x,y)}=yR_2(x^k)$ for some $R_1,R_2\in\KK(x)$. The first factor $\frac{P_1(x,y)}{Q_1(x,y)}$ only depends on $x$, so for $f$ to be birational it must be an element of $\PGLtK$. The conclusion in this case follows.

\emph{Second case.} We now treat the case where $\charK=0$ and where $f$ is of the form $(x,y)\mapsto(\alpha x,y+1)$. Let $(x,y)\dashmapsto(\frac{P_1(x,y)}{Q_1(x,y)},\frac{P_2(x,y)}{Q_2(x,y)})$ be an element of $\Cent(f)$. We have 
\begin{equation}\frac{P_1(\alpha x,y+1)}{Q_1(\alpha x,y+1)}=\frac{\alpha P_1(x,y)}{Q_1(x,y)}\quad \frac{P_2(\alpha x,y+1)}{Q_2(\alpha x,y+1)}=\frac{P_2(x,y)}{Q_2(x,y)}+1.\label{eq:ellipticsecondcase}\end{equation}

The first equation implies that $P_1,Q_1$ are eigenvectors of the $\KK$-linear automorphism $\KK[x,y]\rightarrow\KK[x,y], g(x,y)\mapsto g(\alpha x,y+1)$. We view an element of $\KK[x,y]$ as a polynomial in $x$ with coefficients in $\KK[y]$. Since the only eigenvector of the $\KK$-linear automorphism $\KK[y]\rightarrow\KK[y], g(y)\mapsto g(y+1)$ is $1$ (this is not true if $\charK>0$), we deduce that $P_1,Q_1$ depend only on $x$. Thus, $\frac{P_1(x,y)}{Q_1(x,y)}$ is an element $\eta$ of $\PGLtK$.

We derive $\psi=\frac{P_2}{Q_2}$ and get 
\[\frac{\partial\psi}{\partial y}(\alpha x,y+1)=\frac{\partial \psi}{\partial y}(x,y),\quad \frac{\partial\psi}{\partial x}(\alpha x,y+1)=\alpha^{-1}\frac{\partial \psi}{\partial x}(x,y).\]
As before, this means that $\frac{\partial\psi}{\partial y},\frac{\partial\psi}{\partial x}$ only depend on $x$ (not true if $\charK>0$). Hence, we can write $\psi$ as $ay+B(x)$ with $a\in\KK^*$ and $B\in\KK(x)$. Then equation \eqref{eq:ellipticsecondcase} implies $B(\alpha x)=B(x)+1-a$, which implies further $x\frac{\partial B}{\partial x}(x)$ is invariant under $x\mapsto \alpha x$. If $\alpha$ is of infinite order, then $\frac{\partial B}{\partial x}(x)=\frac{c}{x}$ for some constant $c\in\KK$. This is only possible if $c=0$. So $B$ is constant and $a=1$ in this case. If $\alpha$ is a primitive $k$-th root of unity, then $(\eta(x),ay+B(x))$ commutes with $f^k:(x,y)\mapsto(x,y+k)$. This yields $a=1$ and $B(\alpha x)=B(x)$.

\emph{Third case.} We finally treat the case where $\charK=p>0$ and where $f$ has the form $(x,y)\mapsto(\alpha x,y+1)$ with $\alpha$ of infinite order. Let $g\in\Cent(f)$. Then $g$ commutes with $f^p:(x,y)\mapsto(\alpha^p x,y)$ which is in the form of case 1 (the roles of $x,y$ are exchanged). Thus, we know that $g$ can be written as $(A(y)x,\eta(y))$ where $\eta\in\PGLtK$ and $A\in\KK(x)$. Then $f\circ g=g\circ f$ implies that $\eta$ is $y\mapsto y+R$ for some $R\in\KK$ and that $A(y+1)=A(y)$. The last equation implies $A(y)=S(y)S(y-1)\cdots S(y-p+1)$ for some $S\in\KK(x)$.
\end{pro}

\subsection{\Jonqui twists with trivial action on the base}\label{jonqone}
We follow \cite{CD12} in this section.
\begin{lem}\label{jonqcentjonq}
Let $f\in \Jonq$ be a \Jonqui twist. Let $\Cent(f)$ be the centralizer of $f$ in $\CRK$. Then $\Cent(f)\subset \Jonq$.
\end{lem}
\begin{pro}
The rational fibration preserved by a \Jonqui twist $f$ is unique (cf. \cite{DF01} Lemma 4.5), thus is also preserved by $\Cent(f)$. 
\end{pro}

Let us consider centralizers of \Jonqui twists in $\Jonqz=\operatorname{PGL}_2(\KK(x))$ which is a linear algebraic group over the function field $\KK(x)$. Let $f\in\Jonqz$ and $M=\begin{pmatrix}A&B\\C&D\end{pmatrix}\in \operatorname{GL}_2(\KK(x))$ be a matrix representing $f$ where $A,B,C,D\in \KK[x]$. We introduce the function $\Delta:=\frac{\Trace^2}{\det}$ which is well defined in $\operatorname{PGL}$ and is invariant by conjugation. This invariant $\Delta$ indicates the degree growth:
\begin{lem}[\cite{CD12} Theorem 3.3 \cite{Xie15} Proposition 6.6]
The rational function $\Delta(f)$ is constant if and only if $f$ is an elliptic element.
\end{lem}
\begin{pro}
Let $t_1,t_2$ be the two eigenvalues of the matrix $M$ which are elements of the algebraic closure of $\KK(x)$. The invariant $\Delta(f)$ equals to $t_1/t_2+t_2/t_1+2$. Since $\KK$ is algebraically closed, $\Delta(f)\in\KK$ if and only if $t_1/t_2\in\KK$. If $t_1=t_2$, then by conjugating $M$ to a triangular matrix we can write $f$ in the form $(x,y)\dashmapsto (x,y+a(x))$ with $a\in\KK(x)$ and it follows that $f$ is an elliptic element. 

Suppose now that $t_1\neq t_2$. Let $\zeta:C\rightarrow \PP^1$ be the curve corresponding to the finite field extension $\KK(x)\hookrightarrow \KK(x)(t_1)$, here $\zeta$ is the identity map on $\PP^1$ if $t_1,t_2\in\KK(x)$. The birational transformation $f$ induces a birational transformation $f_C$ on $C\times\PP^1$ by base change. The induced map $f_C$ is of the form $(x,(t_1/t_2)y)$ where $t_1/t_2$ is viewed as a function on $C$. The degree growth of $f_C$ which is the same as $f$ is linear if and only if $t_1/t_2$ is not a constant, i.e.\ if and only if $\Delta(f)$ is not a constant.
\end{pro}

From now on we suppose that $f$ is a \Jonqui twist so that $\Delta(f)\notin \KK$. We still denote by $t_1,t_2$ the two eigenvalues of $M$ as in the above proof, we know that $t_1\neq t_2$. 

We first study the centralizer $\Cent_0(f)$ of $f$ in $\Jonqz=\operatorname{PGL}_2(\KK(x))$. Let $L$ be the finite extension of $\KK(x)$ over which $M$ is diagonalisable; it is $\KK(x)$ itself or a quadratic extension of $\KK(x)$, depending on whether or not $t_1,t_2$ are in $\KK(x)$. The centralizer $\Cent_0^L(f)$ of $f$ in $\operatorname{PGL}_2(L)$ is isomorphic to the multiplicative group $L^*$. So $\Cent_0(f)$, being contained in $\Cent_0^L(f)$ and containing all the iterates of $f$, must be a $1$-dimensional torus over $\KK(x)$. It is split if $L=\KK(x)$, i.e.\ if $t_1,t_2\in\KK(x)$.

If $L=\KK(x)$, then up to conjugation $f$ can be written as $(x,y)\dashmapsto (x,b(x)y)$ with $b\in\KK(x)^*$ and $\Cent_0(f)=\{(x,y)\dashmapsto(x,\gamma(x)y)\vert \gamma\in\KK(x)^*\}$. 

If $L$ is a quadratic extension of $\KK(x)$ and if $\charK\neq 2$, we can put $f$ in a simpler form and write $\Cent_0(f)$ explicitly as follows. We may assume that the matrix $M=\begin{pmatrix}A&B\\C&D\end{pmatrix}$ has entry $C=1$, after conjugation by $\begin{pmatrix}C&0\\0&1\end{pmatrix}$. Once we have $C=1$, a conjugation by $\begin{pmatrix}2&D-A\\0&2\end{pmatrix}$ allows us to put $M$ in the form $\begin{pmatrix}A&B\\1&A\end{pmatrix}$ with $A,B\in\KK[x]$. Therefore $\Cent_0(f)$ is $\{Id,(x,y)\dashmapsto(x,\frac{C(x)y+B(x)}{y+C(x)})\vert C\in\KK(x)\}$ as the ($\KK(x)$-points of the) later algebraic group is easily seen to commute with $f$. Note that $B$ is not a square in $\KK(x)$ because $M$ is not diagonalisable over $\KK(x)$ and that the transformation $f:(x,y)\dashmapsto (x,\frac{A(x)y+B(x)}{y+A(x)})$ fixes pointwise the hyperelliptic curve defined by $y^2=B(x)$.

Now we look at the whole centralizer of $f$. For $\eta\in\PGLtK$ and $f\in\Jonqz$ represented by a matrix $\begin{pmatrix}A(x)&B(x)\\C(x)&D(x)\end{pmatrix}$, we denote by $f_{\eta}$ the element of $\Jonqz$ represented by $\begin{pmatrix}A(\eta(x))&B(\eta(x))\\C(\eta(x))&D(\eta(x))\end{pmatrix}$. Let $f\in\Jonqz$ be a \Jonqui twist and $g:(x,y)\dashmapsto (\eta(x),\frac{a(x)y+b(x)}{c(x)y+d(x)})$ be an element of $\Jonq$. Writing down the commutation equation , we see that $g$ commutes with $f$ if and only if $f$ is conjugate to $f_{\eta}$ in $\operatorname{PGL}_2(\KK(x))$ by the transformation represented by $\begin{pmatrix}a(x)&b(x)\\c(x)&d(x)\end{pmatrix}^{-1}$. We have thus $\Delta(f)(x)=\Delta(f_{\eta})(x)=\Delta(f)(\eta(x))$. Recall that $\Delta(f)\in\KK(x)$ is not in $\KK$. As a consequence the group 
\[
\{\eta\in\PGLtK, \Delta(f)(x)=\Delta(f)(\eta(x))\}
\]
is a finite subgroup of $\PGLtK$. We then obtain:
\begin{thm}[\cite{CD12}]\label{jonqzthm}
Let $f\in\Jonqz$ be a \Jonqui twist preserving the rational fibration fiber by fiber. Let $\Cent(f)$ be the centralizer of $f$ in $\CRK$. Then $\Cent(f)\subset\Jonq$ and $\Cent_0(f)=\Cent(f)\cap\Jonqz$ is a finite index normal subgroup of $\Cent(f)$. The group $\Cent_0(f)$ has a structure of a $1$-dimensional torus over $\KK(x)$. In particular $\Cent(f)$ is virtually abelian.
\end{thm}

\begin{rmk}
In \cite{CD12}, the authors give explicit description of the quotient $\Cent(f)/\Cent_0(f)$ when $\KK=\CC$. 
\end{rmk}

\paragraph{Finite action on the base.} If $f\in\Jonq$ is a \Jonqui twist which has a finite action on the base, then $f^k\in\Jonqz$ for some $k\in\NN$. As $\Cent(f)\subset\Cent(f^k)$, we can use Theorem \ref{jonqzthm} to describe $\Cent(f)$:
\begin{cor}
If $f\in\Jonq$ is a \Jonqui twist which has a finite action on the base, then $\Cent(f)$ is virtually contained in a $1$-dimensional torus over $\KK(x)$. In particular $\Cent(f)$ is virtually abelian. 
\end{cor}
We are contented with this coarse description of $\Cent(f)$ because this causes only a finite index problem as regards the embeddings of $\ZZ^2$ to $\CRK$. We give an example to show how we expect $\Cent(f)$ to look like:
\begin{exa}
If $f$ is $(x,y)\dashmapsto (a(x),R(x)y)$ where $R\in\KK(x)$ and $a\in\PGLtK$ has order $k<+\infty$. Then all maps of the form $(x,y)\dashmapsto (x,S(x)S(a(x))\cdots S(a^{k-1}(x))y)$ with $S\in\KK(x)$ commute with $f$.
\end{exa}

\section{Base-wandering \Jonqui twists}\label{gnljonqtwists}

We introduce some notations. For a Hirzebruch surface $X$, let us denote by $\pi$ the projection of $X$ onto $\PP^1$, i.e.\ the rational fibration. When $X=\PP^1\times\PP^1$, $\pi$ is the projection onto the first factor. For $x\in \PP^1$, we denote by $F_x$ the fiber $\pi^{-1}(x)$. If $f$ is a birational transformation of a Hirzebruch surface $X$ which preserves the rational fibration, we denote by $\overline{f}\in\PGLtK$ the induced action of $f$ on the base $\PP^1$ and we will consider $f$ as an element of $\Jonq$. 

Assume now that $f$ is a \Jonqui twist such that $\overline{f}\in \PGLtK$ if of infinite order, we will call it a \emph{base-wandering \Jonqui twist}. We have an exact sequence:
\begin{equation}
\{1\}\rightarrow \Cent_0(f)\rightarrow \Cent(f) \rightarrow \Cent_b(f)\rightarrow \{1\}
\label{eq:centjonq}
\end{equation}
where $\Cent_0(f)=\Cent(f)\cap \Jonqz$ and $\Cent_b(f)\subset \Cent(\overline{f})\subset \PGLtK$. The action $\overline{f}$ on the base is conjugate to $x\mapsto \alpha x$ with $\alpha\in\KK^*$ of infinite order or to $x\mapsto x+1$. The later case is only possible if $\charK=0$. Thus $\Cent_b(f)$ is a subgroup of $\{x\mapsto \gamma x,\gamma\in\KK^*\}$ or of $\{x\mapsto x+\gamma,\gamma\in\KK\}$. In both cases $\Cent_b(f)$ is abelian. We first remark:
\begin{lem}\label{centzelliptic}
All elements of $\Cent_0(f)$ are elliptic.
\end{lem}
\begin{pro}
By Theorem \ref{jonqzthm}, a \Jonqui twist in $\Jonqz$ cannot have a base-wandering \Jonqui twist in its centralizer. 
\end{pro}
The rest of the article will essentially be occupied by the proof of the following theorem. 
\begin{thm}\label{mainthm}
Let $f\in\Jonq$ be a base-wandering \Jonqui twist. The exact sequence
\[
\{1\}\rightarrow \Cent_0(f)\rightarrow \Cent(f) \rightarrow \Cent_b(f)\rightarrow \{1\}
\]
satisfies
\begin{itemize}
	\item $\Cent_0(f)=\Cent(f)\bigcap\Jonqz$, if not trivial, is $\{(x,y)\mapsto (x,ty),t\in\KK^*\}$, $\{(x,y)\mapsto (x,y+t),t\in\KK\}$ or a torsion group which has order two when $\charK=0$;
	\item $\Cent_b(f)\subset \PGLtK$ is isomorphic to the product of a finite cyclic group with $\ZZ$. The infinite cyclic subgroup generated by $\overline{f}$ has finite index in $\Cent_b(f)$.
\end{itemize}
\end{thm}
\begin{pro}[of Theorem \ref{mainthm}]
The theorem is a consequence of Proposition \ref{centzprop}, Corollary \ref{centpersjonq} and Proposition \ref{centnonpersjonq}. 
\end{pro}

\begin{cor}\label{basewanderingvirtuallyabelian}
When $\charK=0$, the centralizer of a base-wandering \Jonqui twist is virtually abelian.
\end{cor}
\begin{pro}
Theorem \ref{mainthm} implies that $\Cent_b(f)$ is virtually the cyclic group generated by $\overline{f}$ and that $\Cent_0(f)$ is abelian when $\charK=0$. The result follows.
\end{pro}
\begin{rmk}
Theorem \ref{mainthm} is nearly optimal in the sense that $\Cent_b(f)$ can be $\ZZ$ (Remark \ref{exampledeserti}) or a product of $\ZZ$ with a non trivial finite cyclic group (Example \ref{examplecentb}) and $\Cent_0(f)$ can be trivial, isomorphic to $\KK$, $\KK^*$ or $\ZZ/2\ZZ$ (Section \ref{sectioncentz}).
\end{rmk}
\begin{rmk}\label{exampledeserti}
If a base-wandering \Jonqui twist $f$ commutes with an elliptic element of the form $(x,y)\mapsto (x,ty)$ or $(x,y)\mapsto (x,y+t)$ then $f$ can be written as $(\eta(x),yR(x^k))$ or $(\eta(x),y+R(x))$ by Theorem \ref{ellipticthm}. A general base-wandering \Jonqui twist can not be written as $(\eta(x),yR(x^k))$ or $(\eta(x),y+R(x))$. So the centralizer of a general \Jonqui twist $f$ differs from the infinite cyclic group $\left\langle f\right\rangle$ only by some finite groups. For example, for a generic choice of $\alpha,\beta\in \KK^*$, the centralizer of $f_{\alpha,\beta}:(x,y)\dashmapsto(\alpha x,\frac{\beta y+x}{y+1})$ is $\left\langle f_{\alpha\beta}\right\rangle$, this is showed by J. D\'eserti in \cite{Des08}.
\end{rmk}

\subsection{Algebraically stable maps}\label{algebraicallystablemaps}
If $f$ is a birational transformation of a smooth algebraic surface $X$ over $\KK$, we denote by $\Ind(f)$ the set of indeterminacy points of $f$. We say that $f$ is \emph{algebraically stable} if there is no curve $V$ on $X$ such that the strict transform $f^k(V)\subset \Ind(f)$ for some integer $k> 0$. There always exists a birational morphism $\hat{X}\rightarrow X$ which lifts $f$ to an algebraically stable birational transformation of $\hat{X}$ (\cite{DF01} Theorem 0.1). The following theorem says that for $f\in\Jonq$, we can get a more precise algebraically stable model:
\begin{thm}\label{algstabjonq}
Let $f$ be a birational transformation of a ruled surface $X$ that preserves the rational fibration. Then there is a ruled surface $\hat{X}$ and a birational map $\varphi:X\dashmapsto\hat{X}$ such that 
\begin{itemize}
\item the only singular fibers of $\hat{X}$ are of the form $D_0+D_1$ where $D_0,D_1$ are $(-1)$-curves, i.e.\ $\hat{X}$ is a conic bundle;
\item $f_{\hat{X}}=\varphi\circ f\circ\varphi^{-1}$ is an algebraically stable birational transformation of $\hat{X}$ and it preserves the rational fibration of $\hat{X}$ which is induced by that of $X$;
\item $f_{\hat{X}}$ sends singular fibers isomorphically to singular fibers and all indeterminacy points of $f_{\hat{X}}$ and its iterates are located on regular fibers. 
\item $\varphi$ is a sequence of elementary transformations and blow-ups.
\end{itemize}
\end{thm}

Let $z\in X$ be an indeterminacy point of $f$. Let $X\xleftarrow{u}Y\xrightarrow{v}X$ be a minimal resolution of the indeterminacy point $z$, i.e.\ $u,v$ are birational maps which are regular locally outside the fiber over $\pi(z)$, $u^{-1}$ is a series of $n$ blow-ups at $z$ or at its infinitely near points and $n$ is minimal among possible integers.
\begin{lem}\label{jonqindlem}
The total transform by $u^{-1}$ in $Y$ of $F_{\pi(z)}$, the fiber containing $z$, is a chain of $(n+1)$ rational curves $C_0+C_1+\cdots+C_n$: $C_0$ is the strict transform of $F_{\pi(z)}$, $C_0^2=C_n^2=-1$, $C_i^2=-2$ for $0<i<n$ and $C_i\cdot C_{i+1}=1$ for $0\leq i <n$. 
\end{lem}
\begin{pro}
Let us write $u:Y\rightarrow X$ as $Y=Y_n\xrightarrow{u_n}Y_{n-1}\cdots\xrightarrow{u_2}Y_1\xrightarrow{u_1}Y_0=X$ where each $u_i$ is a single contraction of a $(-1)$-curve. Denote by $C_i$ the $(-1)$-curve contracted by $u_i$. By an abuse of notation, we will also use $C_i$ to denote all strict transforms of $C_i$. The connectedness of the fibers and the preservation of the fibration imply that for each $i$, the map $f\circ u_1\circ \cdots \circ u_{i}$ has at most one indeterminacy point on a fiber. To prove the lemma, it suffices to show that the indeterminacy point of $f\circ u_1\circ \cdots \circ u_{i}$ which by construction lies in $C_i$ is not the intersection point of $C_i$ with $C_{i-1}$. 

Suppose by contradiction that $C_{i+1}$ is obtained by blowing up the intersection point of $C_i$ with $C_{i-1}$. Then for $j>i$, the auto-intersection of $C_i$ on $X_j$ is less than or equal to $-2$. Let us write $v:Y\rightarrow X$ as $Y=Y_n\xrightarrow{v_n}Y_{n-1}\cdots\xrightarrow{v_2}Y_1\xrightarrow{v_1}Y_0=X$ where each $v_i$ is a single contraction of a $(-1)$-curve. Since $C_i$ is contracted by $v$, there must exist an integer $k$ such that $v_{k+1}\circ\cdots\circ v_n(C_i)$ is the $(-1)$-curve on $Y_{k}$ contracted by $v_{k}$. By looking at the auto-intersection number, we see that it is possible only if the $C_j,j>i$ are all contracted by $v_{k}\circ\cdots\circ v_n$. But by the minimality of the integer $n$, $C_n$ can not be contracted by $v$.
\end{pro}

\begin{pro}[of Theorem \ref{algstabjonq}]
Our proof is inspired by the proof of Theorem 0.1 of \cite{DF01}.
Let $p_1,\cdots,p_k\in X$ be the indeterminacy points of $f$. 
By Lemma \ref{jonqindlem}, for $1\leq i\leq k$ the minimal resolution of $f$ at $p_i$ writes as
\[
X=X_{i0}\xleftarrow{u_{i1}}X_{i1}\xleftarrow{u_{i2}}\cdots\xleftarrow{u_{in_i}} X_{in_i}=Y_{in_i}\xrightarrow{v_{in_i}}\cdots\xrightarrow{v_{i2}}Y_{i1}\xrightarrow{v_{i1}}Y_{i0}=X
\]
where $u_{i1},\cdots,u_{in_i},v_{i1},\cdots,v_{in_i}$ are single contractions of $(-1)$-curves and $X_{in_i}$ has one singular fiber which is a chain of rational curves $C_{i0}+\cdots+C_{in_i}$. 
Let us write the global minimal resolution of indeterminacy of $f$ by keeping in mind the rational fibration:
\begin{equation*}
\begin{tikzcd}
X=X_0 \arrow[r,dashed]{}{f_0} \arrow{d}{\pi} & X_1 \arrow[r,dashed]{}{f_1} \arrow{d}{\pi} & \cdots \arrow[r,dashed]{}{f_{n-1}} &X_{n} \arrow[r,dashed]{}{f_n} \arrow{d}{\pi} & \cdots \arrow[r,dashed]{}{f_{2n-2}} & X_{2n-1} \arrow[r,dashed]{}{f_{2n-1}} \arrow{d}{\pi} & X_{2n}=X \arrow{d}{\pi} \\
\PP^1 \arrow{r}{\overline{f_0}} & \PP^1 \arrow{r}{\overline{f_1}} & \cdots \arrow{r}{\overline{f_{n-1}}} &\PP^1 \arrow{r}{\overline{f_n}} & \cdots \arrow{r}{\overline{f_{2n-2}}} & \PP^1 \arrow{r}{\overline{f_{2n-1}}} & \PP^1
\end{tikzcd}
\end{equation*}
where $n=n_1+\cdots+n_k$ and
\begin{itemize}
	\item $f_0,\cdots,f_{n-1}$ are blow-ups which correspond to the inverses of $u_{11},\cdots,u_{1n_1},\cdots,u_{k1},\cdots,u_{kn_k}$;
	\item $f_n,\cdots,f_{2n-1}$ are blow-downs which correspond to $v_{11},\cdots,v_{1n_1},\cdots,v_{k1},\cdots,v_{kn_k}$;
	\item $X_n$ has $k$ singular fibers which are chains of rational curves $C_{i0}+\cdots+C_{in_i}, 1\leq i \leq k$;
	\item the abusive notation $\pi$ is self-explaining; we will also denote by $C_{il}$ a strict transform of $C_{il}$ (if it remains a curve) on the surfaces $X_j$. On $X_0=X_{2n}$, it is possible that $C_{i'0}=C_{in_i}$ for $1\leq i,i'\leq k$.
\end{itemize}
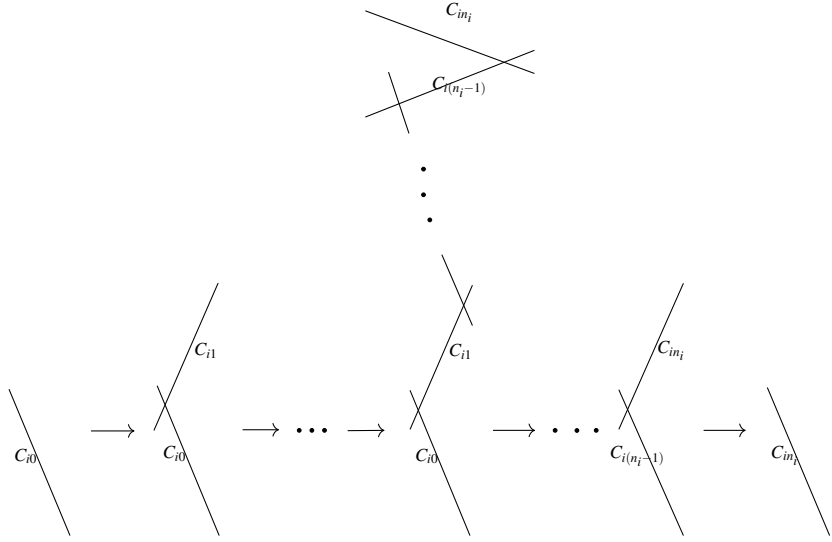
\begin{center}
\begin{tikzpicture}[scale=1.4][line cap=round,line join=round,x=1cm,y=1cm]
\clip(-9.270985411675113,-0.7032400566351622) rectangle (-0.4111856521166486,5.417297941513924);
\draw [line width=0.4pt] (-8.996023788011064,1.3911188248603872)-- (-8.416690318084575,0);
\draw [line width=0.4pt] (-7.589071075332448,1.4242235945704722)-- (-7,0);
\draw [line width=0.4pt] (-7.619907567048753,0.9988180779807584)-- (-7.01,2.4);
\draw [->,line width=0.4pt] (-8.216996768153384,0.9945550742012884) -- (-7.8042473289468015,0.9945550742012884);
\draw [->,line width=0.4pt] (-6.777951426054756,1.005710464450115) -- (-6.4321343283411325,1.005710464450115);
\draw [->,line width=0.4pt] (-5.780116247318101,0.9971761745479083) -- (-5.43,1);
\draw [line width=0.4pt] (-5.1877617757623025,1.3804643620251889)-- (-4.618636891326339,0);
\draw [line width=0.4pt] (-5.1877617757623025,1.008790968107826)-- (-4.595407304206503,2.379336608178102);
\draw [line width=0.4pt] (-4.595407304206503,1.9960484207008213)-- (-4.885777143204444,2.6697064471760417);
\draw [->,line width=0.4pt] (-4.4,1) -- (-4,1);
\draw [line width=0.4pt] (-3.205106242667891,1.386310863707553)-- (-2.5916941560150897,0);
\draw [line width=0.4pt] (-3.205106242667891,1.0096543192716227)-- (-2.5916941560150897,2.3979027259069086);
\draw [->,line width=0.4pt] (-2.4,1) -- (-2,1);
\draw [line width=0.4pt] (-1.795334604921979,1.4078340948181776)-- (-1.203445749379802,0);
\draw [line width=0.4pt] (-5.197875349068657,3.821940607778324)-- (-5.39208472785608,4.404568744140593);
\draw [line width=0.4pt] (-5.612188690481826,3.977308110808262)-- (-4.006724492505796,4.611725414847177);
\draw [line width=0.4pt] (-4.006724492505796,4.391621452221431)-- (-5.612188690481826,4.987196880502862);
\begin{scriptsize}
\draw[color=black] (-8.832275520199392,0.7573428747413243) node {$C_{i0}$};
\draw[color=black] (-7.419843674472681,0.7787433572523351) node {$C_{i0}$};
\draw[color=black] (-7.12023691931853,1.7631655527588315) node {$C_{i1}$};
\draw [fill=black] (-6.24,1) circle (0.5pt);
\draw [fill=black] (-6.12,1) circle (0.5pt);
\draw [fill=black] (-6,1) circle (0.5pt);
\draw[color=black] (-5.022989633239473,0.7573428747413243) node {$C_{i0}$};
\draw[color=black] (-4.701982395574311,1.7524653115033262) node {$C_{i1}$};
\draw [fill=black] (-5,3) circle (0.5pt);
\draw [fill=black] (-5.059999046603208,3.238831331612004) circle (0.5pt);
\draw [fill=black] (-5.059999046603208,3.4827419963702733) circle (0.5pt);
\draw [fill=black] (-3.81,1) circle (0.5pt);
\draw [fill=black] (-3.6032860182144466,0.9988927037163103) circle (0.5pt);
\draw [fill=black] (-3.42,1) circle (0.5pt);

\draw[color=black] (-3.0220445184599645,0.7573428747413243) node {$C_{i(n_i-1)}$};
\draw[color=black] (-2.7010372807948024,1.7631655527588315) node {$C_{in_i}$};
\draw[color=black] (-1.620312913988758,0.7680431159968297) node {$C_{in_i}$};
\draw[color=black] (-4.7126826368298165,4.267022006547094) node {$C_{i(n_i-1)}$};
\draw[color=black] (-4.7126826368298165,4.97323792941045) node {$C_{in_i}$};
\end{scriptsize}
\end{tikzpicture}\end{center}
For any $j\in \NN$, we denote $j\mod 2n$ by $\overline{j}$ and we let $X_j=X_{\overline{j}}$, $f_j=f_{\overline{j}}$.  If $f_j$ blows up a point $r_j\in X_j$, then we denote by $V_{j+1}$ the exceptional curve on $X_{j+1}$. If $f_j$ contracts a curve $W_j\subset X_j$ then we denote by $s_{j+1}$ the point $f_j(W_j)\in X_{j+1}$. For each $V_j$ (resp. $W_j$), there is an $i$ such that $V_j$ (resp. $W_j$) is among $C_{i0},\cdots,C_{in_i}$. 
Suppose that $f$ is not algebraically stable on $H$. Then there exist integers $1\leq M <N$ such that $f_M$ contracts $W_M$ and 
\[
f_{N-1}\circ \cdots \circ f_M(W_M)=r_N\in \Ind(f_N).
\]
Since $f_N$ is a blow-up and $f_M$ is a blow-down, we have $0\leq \overline{N} \leq n-1$ and $n\leq \overline{M} \leq 2n-1$. We can assume that the length $(N-M)$ is minimal. 
Observe first that the minimality of the length implies for all $j$ such that $M\leq j<N-1$, the point $t_{j+1}:=f_j\circ\cdots\circ f_M(W_M)=f_j\circ\cdots\circ f_{M+1}(s_{M+1})$ is neither an indeterminacy point nor a point on a curve contracted by $f_{j+1}$.

The second observation is that for any $j$ such that $M\leq j<N-1$, $t_{j+1}$ is not on a singular fiber of $X_{j+1}$. The reason is as follows. Denote by $M',N'$ the two numbers that satisfy $M< M'\leq N' <N$, $\overline{M'}=\overline{N'}=0$, $M'-M\leq n$ and $N-N' \leq n$. If $M\leq j<M'$ and $t_{j+1}$ was on a singular fiber, then $t_{j+1}$ would be on a component contracted by some $f_k$ with $j+1\leq k<M'$. This contradicts the minimality of $N-M$ by our first observation. Similarly If $N'\leq j<N$ and $t_{j+1}$ was on a singular fiber, then $t_{k}$ would be an indeterminacy point of some $f_k$ with $N'\leq k\leq j$, again contradicting the minimality of $N-M$. If $M'\leq j < N'$ and $t_{j+1}$ was on a singular fiber, then there would exist $L$ such that $M'\leq L\leq j<L+2n\leq N'$ and at least one of the $t_{L},\cdots, t_{L+2n-1}$ is not a point where the corresponding map is locally isomorphic, still contradicting the minimality of $N-M$ by our first observation.

The second observation further implies that for $j$ such that $M\leq j<j+2n <N-1$, $t_{j+1},t_{j+2n+1}$ are not on the same fiber of $X_{j+1}=X_{j+2n+1}$. Suppose to the contrary that this is not true for some $j$. Then there exists $j'$ such that $j<j'\leq j+2n$ and $\overline{j'}=\overline{M}$. $t_{j'}$ would also be on the same fiber as $t_M$. However it is the singular fiber containing $W_M$, contradiction to our second observation that $t_{j'}$ cannot be on a singular fiber.

Since $f_{N-1}$ maps isomorphically the fiber of $X_{N-1}$ containing $t_{N-1}$ (which is regular by the above observation) to the fiber of $X_N$ containing $r_N$, the fiber containing $r_N$ is just one rational curve. As $f_N$ is a blow-up, the fiber of $X_{N+1}$ containing $V_{N+1}$ is the union of two $(-1)$-curves, let us say, $C_{k0}$ and $C_{k1}=V_{N+1}$. Then the fiber of $X_N$ containing $r_N$ is just $C_{k0}$. Similarly the singular fiber of $X_M$ containing $W_M$ is $C_{mn_m}+C_{m(n_m-1)}$ for some $1\leq m\leq k$ and $W_M=C_{m(n_m-1)}$. 

\emph{First case.} Suppose that $m=k$ and $n_k=1$. Let $a\in\NN$ be the minimal integer such that $M+2an>N$. Then for $N<j\leq M+2an$, the surface $X_j$ has a singular fiber $C_{k0}+C_{k1}$ and the maps $f_N,\cdots,f_{M+2an-1}$ are all regular on $C_{k0}+C_{k1}$.
Now we blow-up $t_{M+1},\cdots,t_{N-1},r_N$. For $\overline{j_1}=\overline{j_2}$, we showed that $t_{j_1},t_{j_2}$ are not on the same fiber of $X_{j_1}=X_{j_2}$. This means that these blow-ups only give rise to singular fibers which are unions of two $(-1)$-curves. We denote by $\hat{X}_j$ the modified surfaces, and $\hat{f}_j$ the induced maps. Then every $\hat{X}_j$ has singular fibers of the form $C_{k0}+C_{k1}$ and every $\hat{f}_j$ is regular around these singular fibers. Let $\hat{f}=\hat{f}_{2n-1}\circ \cdots \circ \hat{f}_0$. The number of indeterminacy points of $\hat{f}$ (it was $k$ for $f$) has decreased by one. Note that $\hat{f}$ exchanges the two components $C_{k0}$ and $C_{k1}$. This fact will be used in the proof of Proposition \ref{finiteorderjonq}.

\emph{Second case.} Suppose that $m=k$ and $n_k>1$ or simply $m\neq k$.
We blow-up $r_N$ and contract the strict transform of the initial fiber containing $r_N$ which is $C_{k0}$, obtaining a new surface $\hat{X}_N$ whose corresponding fiber is now the single rational curve $C_{k1}$. 
We perform elementary transformations at $t_{N-1},\cdots,t_{M+1}$, i.e.\ we blow-up $X_j$ at $t_j$ and contract the strict transform of the initial fiber, replacing $X_j$ with $\hat{X}_j$. This process has no ambiguity: if $\overline{j_1}=\overline{j_2}$, we showed that $t_{j_1},t_{j_2}$ are not on the same fiber of $X_{j_1}=X_{j_2}$, so the corresponding elementary transformations do not interfere with each other. Let us denote by $\hat{f}_M,\cdots,\hat{f}_N$ the maps induced by $f_M,\cdots,f_N$.  

We now analyse the effects of $\hat{f}_M,\cdots,\hat{f}_N$. First look at $f_N$, it lifts to a regular isomorphism after blowing up $r_N$. Thus $\hat{f}_N$ is the blow-up at the point $e_N$ of $\hat{X}_N$ to which $C_{k0}$ is contracted. After this step, the map going from $X_{N-1}$ to $\hat{X}_N$ induced by $f_{N-1}$ is as following: it contracts the fiber containing $t_{N-1}$ to $e_N$ and blows up $t_{N-1}$. Then we make elementary transformations at $t_{N-1},\cdots,t_{M+1}$ in turn. The maps $\hat{f}_{N-1},\cdots,\hat{f}_{M+1}$ are all regular on the modified fibers, thus they are still single blow-ups or single blow-downs. The behaviour of $\hat{f}_{M}$ differs from the previous ones: it does not contract $C_{m(n_m-1)}$ any more, but contracts $C_{mn_m}$. 

The hypothesis $m\neq k$ (or $m=k$, $n_k>1$) forbids $C_{k0}\subset X_{N+1}$ to go back into the fiber of $X_{M+2na}=X_M$ containing $W_M$ without being contracted. More precisely this implies the existence of $N'>N$ such that 
\begin{itemize}
\item $X_{N+1},\cdots,X_{N'}$ all contain $C_{k0}$ and $C_{k1}$;
\item $f_{N+1},\cdots,f_{N'-1}$ are regular on $C_{k0}$ and $f_{N'}$ contracts $C_{k0}$;
\item if $a\in \NN$ is the minimal integer such that $M+2na>N$, then $N'<M+2na$.
\end{itemize}
On the surfaces $X_{N+1},\cdots,X_{N'}$, $C_{k0}$ is always a $(-1)$-curve, we contract all these $C_{k0}$ and obtain new surfaces $\hat{X}_{N+1},\cdots,\hat{X}_{N'}$. The second and the third property listed above mean that the new induced maps $\hat{f}_{N},\cdots,\hat{f}_{N'}$ are all single blow-ups, single blow-downs or simply isomorphisms. 

In summary we get a commutative diagram:
\begin{equation*}
\begin{tikzcd}
\hat{X}_0 \arrow[r,dashed]{}{\hat{f}_0} \arrow[d,dashed]{} & \hat{X}_1 \arrow[r,dashed]{}{\hat{f}_1} \arrow[d,dashed]{} & \cdots \arrow[r,dashed]{}{\hat{f}_{n-1}} &\hat{X}_{n} \arrow[r,dashed]{}{\hat{f}_n} \arrow[d,dashed]{} & \cdots \arrow[r,dashed]{}{\hat{f}_{2n-2}} & \hat{X}_{2n-1} \arrow[r,dashed]{}{\hat{f}_{2n-1}} \arrow[d,dashed]{} & \hat{X}_{2n}=\hat{X}_0 \arrow[d,dashed]{} \\
X_0 \arrow[r,dashed]{}{f_0} & X_1 \arrow[r,dashed]{}{f_1} & \cdots \arrow[r,dashed]{}{f_{n-1}} &X_{n} \arrow[r,dashed]{}{f_n} & \cdots \arrow[r,dashed]{}{f_{2n-2}} & X_{2n-1} \arrow[r,dashed]{}{f_{2n-1}} & X_{2n}=X_0\end{tikzcd}
\end{equation*}
where the vertical arrows are composition of elementary transformations and blow-ups. 
Let us remark that:
\begin{itemize}
  \item the first vertical arrow $\hat{X}_0\dashrightarrow X_0$ is a composition of elementary transformations.
	\item the blow-ups or the contractions of the $\hat{f}_j$ only concern the $k$ singular fibers and the exceptional curves are always among $C_{10},\cdots,C_{1n_1},\cdots,C_{k0},\cdots,C_{kn_k}$;
	\item there is no more $C_{k0}$. We then do a renumbering: $C_{k1},\cdots,C_{kn_k}$ become $C_{k0},\cdots,C_{k(n_k-1)}$.
\end{itemize} 

Let $\hat{f}=\hat{f}_{2n-1}\circ \cdots \circ \hat{f}_0$. We repeat the above process. Either we are in the first case and $k$ decreases, or we are in the second case and the total number of $C_{10},\cdots,C_{1n_1},\cdots,C_{k0},\cdots,C_{kn_k}$ decreases. As a consequence, after a finite number of times, either we get an algebraically stable map $\hat{f}$, or we will get rid of all the $C_{10},\cdots,C_{1n_1},\cdots,C_{k0},\cdots,C_{kn_k}$. In the later case $\hat{f}$ is a regular automorphism, thus automatically algebraically stable. 
\end{pro}

\subsection{Elliptic elements in the \Jonqui group}
Some materials in this section are taken from \cite{Bla11}. We assume $\charK=0$ in this section.

\begin{lem}\label{telescopic}
Assume $\charK=0$.
Let $f:(x,y)\dashmapsto(\eta(x),yR(x)), \eta\in\PGLtK, R\in\KK(x)$ be an elliptic element. Then
\begin{enumerate}
	\item either $R\in\KK$,
	\item or $R(x)=\frac{rS(x)}{S(\eta(x))}$ with $r\in\KK^*$ and $S\in\KK(x)\backslash \KK$.
\end{enumerate}
\end{lem}
\begin{pro}
If $\eta$ is the identity, then we see easily, by looking at the degree growth, that $f$ is elliptic if and only if $R$ is constant.
From now on assume that $\eta$ is not the identity. 

We first consider the case where $f$ has infinite order. Then by Proposition \ref{ellipticnormalform} $f$ is conjugate to an automorphism of a Hirzebruch surface. By Theorem \ref{algstabjonq}, the conjugation which turns $f$ into an automorphism of a Hirzebruch surface is a sequence of elementary transformations. After conjugation it preserves the two strict transforms of the two sections $\{y=0\}$ and $\{y=\infty\}$.
Therefore there exists $g\in \Jonqz$ of the form $(x,y)\dashmapsto(x,yS(x)), S\in\KK(x)$ such that $g\circ f\circ g^{-1}$ is $(x,y)\dashmapsto(\eta(x),r y)$ with $r\in \KK^*$. Hence $f$ is $(x,y)\dashmapsto(\eta(x),y\frac{rS(x)}{S(\eta(x))}$. 

Now assume that $f$ has finite order. Then $\eta$ has also finite order. Denote by $d$ the order of $\eta$ and $dk$ the order of $f$. The expression of $f^d$ is $(x,y)\dashmapsto(x,T(x) y)$ where 
\[T(x)=R(x)R(\eta(x))\cdots R(\eta^{d-1}(x)).\]
As $f^{dk}=\Id$, we have $T(x)^k=1$. Therefore $T(x)=\alpha$ is a root of unity. Take a $\beta\in \KK^*$ such that $\beta^d=\alpha$. Then $R'(x)=R(x)/\beta$ satisfies 
\[R'(x)R'(\eta(x))\cdots R'(\eta^{d-1}(x))=1.\]
In other words $R'$ represents a cocycle of $\ZZ/d\ZZ$ with values in $\KK^*(X)$. By Hilbert's Theorem 90, the cohomology group $\operatorname{H}^1(\ZZ/d\ZZ,\KK^*(X))$ is trivial. Thus the existence of $S\in \KK^*(X)$ such that $R'(x)=\frac{S(x)}{S(\eta(x))}$. The conclusion follows.
\end{pro}

\begin{rmk}
If $\charK>0$ then we cannot apply Hilbert's Theorem 90 in the above situation because the corresponding field extension may not be seperable.
\end{rmk}

\begin{rmk}
In the above lemma $S$ may not be unique. If $\eta$ has finite order and $T\in\KK(x)$ is such that $T(x)=T(\eta(x))$, then $\frac{S(x)}{S(\eta(x))}=\frac{T(x)S(x)}{T(\eta(x))S(\eta(x))}$. 
\end{rmk}
A direct corollary of Lemma \ref{telescopic} is:
\begin{cor}\label{diagonalisable}
If $\charK=0$ then a diagonalisable elliptic element of $\Jonqz=\operatorname{PGL}_2(\KK(X))$ is always conjugate to an automorphism of a Hirzebruch surface.
\end{cor}

\begin{lem}\label{Risapolynomial}
Let $f:(x,y)\dashmapsto(\eta(x),y+R(x)), \eta\in\PGLtK, R\in\KK(x)$ be an elliptic element. Then one of the following holds:
\begin{enumerate}
	\item $\eta$ has finite order;
	\item there exists a coordinate $x'$ such that $\eta$ is $x'\mapsto \nu x'$ with $\nu\in\KK^*$, and $R(x')=R_1(x')+R_2(\frac{1}{x'})$ where $R_1,R_2$ are polynomials in $x'$.
	\item there exists a coordinate $x'$ such that $\eta$ is $x'\mapsto x'+1$, and $R$ is a polynomial in $x'$.
\end{enumerate}
\end{lem}
\begin{pro}
Assume that $\eta$ has infinite order, then for some coordinate $x'$, $\eta$ can be written as $\eta(x')=x'\mapsto x'+1$ or $x'\mapsto \nu x'$ with $\nu\in\KK^*$. In coordinates $(x',y)$, write the transformation $f$ as $(x',y)\dashmapsto(\eta(x'),y+C(x')+\frac{P(x')}{Q(x')})$ where $C,P,Q\in \KK[x']$ are such that either $P=0$ or $\deg P<\deg Q$. There is nothing to prove if $P=0$. Suppose by contradiction that $\deg P<\deg Q$. For $n\in\NN^*$, the iterate $f^n$ is 
\[
(x',y)\dashmapsto\left(\eta^n(x'),y+\Sigma(x')+\frac{P(x')}{Q(x')}+\cdots+\frac{P(\eta'^{n-1}(x'))}{Q(\eta'^{n-1}(x'))}\right)
\] 
Where $\Sigma(x')$ is a polynomial of degree at most the degree of $C$.
If $\eta(x')=\nu x'$ and $Q(x')=qx'^d$ for some $q\in \KK$ and $d\in\NN$, then replacing the coordinate $x'$ with $x^*=1/x'$ we can write $f$ as $(x^*,y)\dashmapsto(\frac{x^*}{\nu},y+C(\frac{1}{x^*})+q^{-1}x^{*d}P(1/x^*))$. As $d=\deg Q>\deg P$, the coordinate $x^*$ satisfies the required properties.

In the remaining cases, either $\eta(x')=x+1$ or $Q$ has a factor $x-a$ with $a\neq 0$. In both cases some factor of $Q$ gives rise to an infinite number of distinct factors in the polynomials $Q(\eta'^{n}(x')), n\in \NN$. Then if we write 
\[\frac{P(x')}{Q(x')}+\cdots+\frac{P(\eta'^{n-1}(x'))}{Q(\eta'^{n-1}(x'))}\] as a single fraction, then the degree of its denominator tends to infinity as $n$ tends to infinity. This contradicts the hypothesis that the degrees of the $f^n$s are bounded. This means $P=0$ in these situations.
\end{pro}

We refer to \cite{Bla11} Section 3 (especially Proposition 3.3 and Lemma 3.9) for results that are stronger than the following proposition. 
\begin{prop}\label{finiteorderjonq}
Assume that $\charK=0$.
Let $f\in\Jonq$ be an elliptic element that is not conjugate to an automorphism of a Hirzebruch surface. Then $f$ has even finite order $2k$ and is conjugate to an automorphism $f'$ of a conic bundle. Moreover $f^k$ is in $\Jonqz$ and ${f'}^k$ exchanges the two components of at least one singular fiber of the conic bundle. 
\end{prop}
\begin{pro}
We see by Proposition \ref{ellipticnormalform} that an elliptic element of infinite order is always conjugate to an automorphism of Hirzebruch surface. Hence our hypothesis implies immediately that $f$ is of finite order. We can assume that $f'$ is an algebraically stable map on a conic bundle $X$ which satisfies the conclusions of Theorem \ref{algstabjonq}. We claim that $f'$ is an automorphism of $X$. Suppose by contradiction that $p_0$ is an indeterminacy point of $f'$. It must lie on a regular fiber $F$ of $X$. Let $d$ be the order of $f'$. As $f'$ is algebraicaly stable, $p_0$ is not an indeterminacy point of ${f'}^{-n}$ for any $n$ and there exists $p_1, \cdots, p_d$ belonging to the fibers ${f'}^{-1}(F),\cdots {f'}^{-d}(F)$ such that $f'$ is well defined at each $p_j$ and $f'(p_j)=p_{j-1}$. However $f'^d=\Id$ implies $p_d=p_0$, contradicting that $p_0$ is an indeterminacy point.

Since by hypothesis $X$ is not a Hirzebruch surface, it must have some singular fibers. As $f'$ has finite order, there exists a minimal $k>0$ such that ${f'}^k$ preserves the rational fibration fiber by fiber, i.e.\ $f^k\in \Jonqz$. 

Let $M=\begin{pmatrix}A(x)&B(X)\\C(X)&D(X)\end{pmatrix}$ be a matrix of $\Gl_2(\KK(X))$ that represents $f^k\in\Jonqz=\operatorname{PGL}_2(\KK(X))$. 

We claim that $M$ is not diagonalisable. It is proved when $\KK=\CC$ in \cite{Bla11}[page 479 parts (i), (ii) of the proof of Proposition 3.3] that this would imply that $f$ is conjugate to an automorphism of a Hirzebruch surface; the proof works word by word when $\charK=0$ . A consequence is that ${f'}^k$ exchanges the two components of at least one singular fiber of $X$ because otherwise ${f'}^k$ would be conjugate to an automorphism of a Hirzebruch surface and would be diagonalisable (see the formula \eqref{eq:autofhirz}). Thus the order of $f^k$ is even.

Let $2j$ be the order of $f^k$. Then there exists $\lambda\in \KK(x)^*$ such that $M^{2j}=\lambda \Id$ is a scalar matrix. Taking determinants of both side we get $\det M^{2j}=\lambda^2$. Thus there exists $\kappa\in \KK(X)$ such that $\kappa^j=\lambda$. Let $P_M\in \KK(X)[T]$ be the minimal polynomial of $M$. Then $P_M$ is a factor of $T^{2j}-\lambda=T^{2j}-\kappa^j$. Since $M$ is not diagonalisable the polynomial $P_M$ is irreducible in $\KK(X)[T]$. Any irreducible factor of $T^{2j}-\kappa^j$ has the form $T^2-\mu \kappa$ for some $j$-th root of unity $\mu$. Therefore $M^2$ is a scalar matrix. In other words $j=1$ and $f^k$ is an involution.
\end{pro}

\begin{cor}\label{ordertwoboy}
Assume that $\charK=0$.
Let $f\in\Jonqz$ be an elliptic element which is not conjugate to an automorphism of a Hirzebruch surface. Then $f$ is an involution and is conjugate to an automorphism $f'$ of a conic bundle $X$. The automorphism $f'$ fixes pointwise a hyperelliptic curve whose projection onto the base $\PP^1$ is a ramified double cover. In some affine chart $f$ can be written as $(x,y)\dashmapsto (x,\frac{a(x)}{y})$ with $a\in\KK[x]$. The hyperelliptic curve is given by the equation $y^2=a(x)$.
\end{cor}
Such involutions are well known and are called \emph{\Jonqui involutions}, see \cite{BB00}.

\subsection{The group $\Cent_0(f)$}\label{sectioncentz}
Now we turn back to the study of centralizers.
Let $f$ be a base-wandering \Jonqui twist. In \cite{CD12}, it is proved by explicit calculations, in the case where $\KK=\CC$, that $\Cent_0(f)$ is isomorphic to $\CC^*$, $\CC^*\rtimes\ZZ/2\ZZ$, $\CC$ or a finite group (this is not optimal). Their arguments do not work directly when $\charK>0$. With a more precise description of elements of $\Jonqz$, we simplify their arguments and improve their results.

Let $g\in\Cent_0(f)$ be non trivial. By Lemma \ref{centzelliptic} $g$ is elliptic. If $g$ is conjugate to an automorphism of a Hirzebruch surface, then by proposition \ref{jonqzelliptic} by choosing suitable coordinates we can write $g$ as $(x,y)\mapsto (x,\beta y)$ or $(x,y)\mapsto (x,y+1)$. 

\begin{lem}
Suppose that there exists a non trivial $g\in\Cent_0(f)$ that can be written as $(x,y)\mapsto (x,\beta y)$ with $\beta\in\KK^*$. Either $f$ has the form $(a(x),R(x)y^{-1})$ with $a\in \operatorname{PGL}_2(\KK)$ and $\Cent_0(f)$ is an order two group generated by the involution $(x,y)\mapsto(x,-y)$, or $f$ has the form $(a(x),R(x)y)$ and $\Cent_0(f)$ is $\{(x,y)\mapsto (x,\gamma y),\gamma\in\KK^*\}$.
\end{lem}
\begin{pro}
The map $g$ preserves $\{y=0\}$ and $\{y=\infty\}$ and these two curves are the only $g$-invariant sections. Thus $f$ permutes these two sections and has necessarily the form $(x,y)\dashmapsto(a(x),R(x)y^{\pm 1})$ where $R\in\KK(x)$ and $a\in\PGLtK$ has infinite order. If $f$ is $(a(x),R(x)y^{-1})$, then $\beta=-1$. Since $\Cent(f)\subset \Cent(f^2)$ we replace $f$ by $f^2$ in what follows so that we can assume $f$ is $(a(x),R(x)y)$. 

The only $f$-invariant sections are $\{y=0\}$ and $\{y=\infty\}$. Indeed an invariant section $s$ satisfies 
\[
s(a^n(x))=R(x)\cdots R(a^{n-1}(x))s(x) \quad \forall n\in\NN.
\]
If $s$ was not $\{y=0\}$ nor $\{y=\infty\}$, then the two sides of the equations would be non-zero rational fractions and the degree of the right side would go to infinity with $n$ because $f$ is a \Jonqui twist. Thus, an element of $\Cent_0(f)$ permutes the two $f$-invariant sections and has the form $(x,A(x)y)$ or $(x,\frac{A(x)}{y})$ with $A\in\KK(x)$. In the first case the commutation relation implies $A(a(x))=A(x)$ which further implies that $A$ is a constant. In the second case the commutation relation gives $A(a(x))^{-1}R(x)^2A(x)=1$ which further implies that $(a(x),R(x)^2 y)$ is conjugate by $(x,A(x)y)$ to an elliptic element $(a(x), y)$. This is not possible because the map $f':(x,y)\dashmapsto(a(x),R(x)^2 y)$ is a \Jonqui twist. Indeed the iterates $f^n,f'^n$ are respectively 
\[
(a^n(x), R(x)\cdots R(a^{n-1}(x))y)\quad \text{and} \quad (a^n (x), (R(x)\cdots R(a^{n-1}(x)))^2 y)
\] 
and they have the same degree growth.

Reciprocally all elements of the form $(x,y)\mapsto (x,\beta y)$ with $\beta\in\KK^*$ commute with $f:(x,y)\dashmapsto(a(x),R(x)y)$ and we have already observed that $(x,y)\mapsto(x,-y)$ is the only non trivial element of $\Jonqz$ which commutes with $(a(x),R(x)y^{-1})$.
\end{pro}

\begin{lem}
Suppose that there exists a non trivial $g\in\Cent_0(f)$ that can be written as $(x,y)\mapsto (x,y+1)$. Then $f$ has the form $(a(x),y+S(x))$ with $S\in\KK(x)$ and $\Cent_0(f)$ is $\{(x,y+\gamma),\gamma\in\KK\}$.
\end{lem}
\begin{pro}
The section $\{y=\infty\}$ is the only $g$-invariant section. Thus $f$ preserves this section and has the form $(x,y)\dashmapsto(a(x),R(x)y+S(x))$ where $R,S\in\KK(x)$ and $a\in\PGLtK$ is of infinite order. Writing down the relation $f\circ g=g\circ f$, we see that $R=1$. Thus $f$ is $(a(x),y+S(x))$ where $S$ belongs to $\KK(x)$ but not to $\KK[x]$ since $f$ is a \Jonqui twist. The only $f$-invariant section is $\{y=\infty\}$. Indeed an invariant section $s$ satisfies 
\[
s(a^n(x))=s(x)+S(x)+\cdots +S(a^{n-1}(x)) \quad \forall n\in\NN.
\]
If $s$ was not $\{y=\infty\}$, then the two sides of the equations are rational fractions. The degree of the right-hand side grows linearly in $n$ while the degree of the left-hand side does not depend on $n$, contradiction. Thus, an element of $\Cent_0(f)$ fixes $\{y=\infty\}$ and has the form $(x,A(x)y+B(x))$ with $A,B\in\KK(x)$. Writing down the commutation relation, we get 
\[
A(x)y+B(x)+S(x)=A(a(x))y+A(a(x))S(x)+B(a(x)).
\] 
The fact that $a$ has infinite order implies that $A$ is a constant. Then the equation is reduced to 
\[
B(x)+(1-A)S(x)-B(a(x))=0.
\]
If $A\neq 1$, then $f:(x,y)\dashmapsto(a(x),y+S(x))$ would be conjugate by $(x,y+\frac{B(x)}{1-A})$ to the elliptic elment $(a(x),y)$. Therefore $A=1$ and $B$ is a constant. 
Reciprocally we see that all elements of the form $(x,y)\mapsto (x,y+\beta)$ with $\beta\in\KK$ commute with $f:(x,y)\dashmapsto(a(x),y+S(x))$.
\end{pro}

\begin{lem}
Assume that no non-trivial element of $\Cent_0(f)$ is conjugate to an automorphism of a Hirzebruch surface and that $\Cent_0(f)$ has a non-trivial element $g$. Then $\Cent_0(f)$ is a torsion group. If moreover $\charK=0$ then $g$ has order two and is the only non-trivial element of $\Cent_0(f)$.
\end{lem}
\begin{pro}
By Lemma \ref{centzelliptic}, $g$ is an elliptic element. Since an elliptic element of infinite order is conjugate to an automorphism of a Hirzebruch surface by Proposition \ref{ellipticnormalform}, $g$ has finite order. If $\charK=0$ then $g$ is a \Jonqui involution of a conic bundle $X$ by Corollary \ref{ordertwoboy}.

Assume now that $\charK=0$. Then by Corollary \ref{ordertwoboy} $g$ fixes pointwise a hyperelliptic curve $C$. The map $f$ induces an action on $C$, equivariant with respect to the ramified double cover. The action of $f$ on $C$ is infinite, this is possible only if the action of $f$ on the base is up to conjugation $x\mapsto \alpha x$ and if $C$ is a rational curve whose projection on the base $\PP^1$ is ramified over $x=0,x=\infty$. Then the only singular fibers of $X$ are over $x=0,x=\infty$. Note that $f$ can have at most one indeterminacy point on each of these two fibers. If $f$ had an indeterminacy point on these two fibers, then it would be a fixed point of $g$ because $g$ commutes with $f$ and $g$ preserves each fiber. But the only fixed point of $g$ on a singular fiber is the intersection point of the two components, which can not be an indeterminacy point by Lemma \ref{jonqindlem}. Therefore the \Jonqui twist $f$ must have an indeterminacy point over a point whose orbit in the base is infinite. This implies that the indeterminacy points of all the iterates of $f$ form an infinite set. As $g$ commutes with all the iterates of $f$, it fixes an infinite number of these indeterminacy points. Thus, the hyperelliptic curve $C$ associated to $g$ is the Zariski closure of these indeterminacy points and is uniquely determined by $f$. However $C$ determines $g$ too. (cf. \cite{BB00} Proposition 2.7). Therefore $g$ is uniquely determined by $f$ and is the only non trivial element of $\Cent_0(f)$. 
\end{pro}

Putting together the three previous lemmas, we obtain the following improvement of \cite{CD12}:
\begin{prop}\label{centzprop}
Let $f$ be a base-wandering \Jonqui twist. If $\Cent_0(f)$ is not trivial, then it is $\{(x,y)\mapsto (x,ty),t\in\KK^*\}$, $\{(x,y)\mapsto (x,y+t),t\in\KK\}$, $\langle (x,y)\mapsto (x,-y)\rangle$ or a torsion group. If $\charK=0$ and $\Cent_0(f)$ is a non-trivial torsion group, then it is a group of order two generated by a \Jonqui involution.
\end{prop}

\subsection{Persistent indeterminacy points}
\subsubsection{general facts}
Let $f$ be a birational transformation of a surface $X$. An indeterminacy point $x\in X$ of $f$ will be called \emph{persistent} if 1) for every $i>0$, $f^{-i}$ is regular at $x$; and 2) there are infinitely many curves contracted onto $x$ by the iterates $f^{-n}$, $n\in \NN$. This notion of persistence and the following idea appeared first in a non published version of \cite{Can11}, and it is also applied to some particular examples in \cite{Des08}. Note that it is different from the definition of persistent base points given in \cite{BD15}.
\begin{prop}\label{orbitargument}
Let $f$ be an algebraically stable birational transformation of a surface $X$. Suppose that there exists at least one persistent indeterminacy point with an infinite backward orbit. Let $n$ denote the number of such indeterminacy points. Then the centralizer $\Cent(f)$ of $f$ admits a morphism $\varphi:\Cent(f)\rightarrow \Sym_n$ to the symmetric group of order $n$ satisfying the following property: for any $g\in\Ker(\varphi)$, there exists $l\in\ZZ$ such that $g\circ f^l$ preserves fiber by fiber a pencil of rational curves.
\end{prop}
\begin{pro}
The algebraic stability of $f$ will be used throughout the proof, we will not recall it each time.
Denote by $p_1,\cdots,p_n$ the persistent indeterminacy points of $f$. Let $g$ be a birational transformation of $X$ which commutes with $f$. Fix an index $1\leq n_0\leq n$. Since $\{f^{-i}(p_{n_0}),i>0\}$ is infinite, there exists $k_0>0$ such that $g$ is regular at $f^{-k}(p_{n_0})$ for all $k\geq k_0$. For infinitely many $j>0$, $f^{-j}$ contracts a curve onto $p_{n_0}$, denote these curves by $C_{n_0}^j$. There exists $k_1>0$ such that $g$ does not contract $C_{n_0}^k$ for all $k\geq k_1$. We deduce, from the above observations and the fact that $f$ and $g$ commute, that for $k\geq k_0$ the point $g(f^{-k}(p_{n_0}))$ is an indeterminacy point of some $f^m$ with $0<m\leq k+k_1$. Then there exists $0\leq m_0< m$ such that
\begin{itemize}
	\item for $0\leq i\leq m_0$, $f^i$ is regular at $g(f^{-k}(p_{n_0}))$;
	\item $f^{m_0}(g(f^{-k}(p_{n_0})))=g(f^{m_0-k}(p_{n_0}))$ is an indeterminacy point of $f$.
\end{itemize}
By looking at $g(f^{-k}(p_{n_0}))$ and $C_{n_0}^{k'}$ for infinitely many $k,k'$, we see that the above indeterminacy point does not depend on $k$ and is persistent with an infinite backward orbit. So it is $p_{\sigma_g(n_0)}$ for some $1\leq \sigma_g(n_0)\leq n$. 
This gives us a well defined map $\sigma_g:\{1,\cdots,n\}\rightarrow \{1,\cdots,n\}$.

Now let $g,h$ be two elements of $\Cent(f)$. Then by considering a sufficiently large $k$ for which $g$ is regular at $f^{-k}(p_{n_0})$ and $h$ is regular at $g(f^{-k}(p_{n_0}))$, we see that $\sigma_h\circ\sigma_g=\sigma_{h\circ g}$. By taking $h=g^{-1}$ we see that $\sigma_g$ is bijective. We have then a group homomorphism $\varphi$ from $\Cent(f)$ to the symmetric group $\Sym_n$ which sends $g$ to $\sigma_g$.

Assume that $n_0$ is a fixed point of $\sigma_g$, this holds in particular when $g\in \Ker(\varphi)$. We keep the previous notations. Since $g(f^{-k}(p_{n_0}))$ is an indeterminacy point of $f^m$ whose forward orbit meets $p_{n_0}$, for an appropriate choice of $l\leq k$ we have
\[
	g\circ f^l(f^{-k}(p_{n_0}))=f^{-k}(p_{n_0}) 
\]
for all $k\geq k_0$. This implies further
\[
g\circ f^l(C_{n_0}^{k'})=C_{n_0}^{k'}  
\]
for all sufficiently large $k'$. We conclude by Lemma \ref{invhypersurfaces} below.
\end{pro}

The proof of the following lemma in \cite{Can10} is written over $\CC$ for rational self-maps. The same proof works in all characteristics for birational self-maps but not for general rational self-maps (this is also observed and used in \cite{Xie15}.). 
\begin{lem}\label{invhypersurfaces}
If a birational transformation $f$ of a smooth algebraic surface preserves infinitely many curves then these curves are members of a pencil of curves and $f$ preserves each member of this pencil of curves.
\end{lem}
Note that Lemma \ref{invhypersurfaces} is not needed in the proof of Proposition \ref{orbitargument} if $f\in \Jonq$, and this is the only case in this paper where we apply Proposition \ref{orbitargument} (see Corollary \ref{centpersjonq}).

\begin{rmk}
Using the tools introduced in \cite{LU} and the same idea of Proposition \ref{orbitargument}, Lonjou-Urech generalized Proposition \ref{orbitargument} and Corollary \ref{centpersjonq} below to higher dimension, see \cite{LU} Theorem 1.6. 
\end{rmk}

\subsubsection{persistent indeterminacy points for \Jonqui twists}
We examine the notion of persistence in the \Jonqui group and give a complement to Theorem \ref{algstabjonq}:
\begin{prop}\label{persjonq}
Let $f$ be a \Jonqui twist acting algebraically stably on a conic bundle $X$ as in the statement of Theorem \ref{algstabjonq}. Then an indeterminacy point $p$ of $f$ is persistent if and only if the orbit of $\pi(p)\in \PP^1$ under $\overline{f}$ is infinite. And in that case, every $f^{-i}, i\in \NN^*$ contracts a curve onto $p$.
\end{prop}
\begin{pro}
If $\pi(p)$ has a finite orbit then $p$ certainly cannot be persistent. Let us assume that the orbit of $\pi(p)$ is infinite. Then $\overline{f}$ is conjugate to $x\mapsto \alpha X$ with $\alpha\in K^*$ of infinite order or to $x\mapsto x+1$ (only when $\charK=0$). By the algebraic stability of $f$, $f^{-i}$ is regular at $p$ for all $i>0$ and all the points $f^{-i}(p),i>0$ are on distinct fibers. Denote by $x_0,x_1$ the points $\pi(p),\overline{f}(\pi(p))$. By Theorem \ref{algstabjonq}, we know that the fibers $F_{x_0},F_{x_1}$ are not singular. Thus $f$ is regular on $F_{x_0}\backslash \{p\}$ and contracts it onto a point $q\in F_{x_1}$; $f^{-1}$ is regular on $F_{x_1}\backslash \{q\}$ and contracts it onto $p$. Now pick a point $x_n$ in the forward orbit of $x_0$ by $\overline{f}$ and consider the fiber $F_{x_n}$. The fiber $F_{x_n}$ cannot be contracted onto $q$ by $f^{-(n-1)}$ because of the algebraic stability of $f$. As a consequence it is contracted by $f^{-n}$ onto $p$. 
\end{pro} 

\begin{cor}\label{centpersjonq}
Let $f$ be a \Jonqui twist acting algebraically stably on a conic bundle $X$ as in the statement of Theorem \ref{algstabjonq}. Suppose that the base action $\overline{f}\in\PGLtK$ has infinite order and there is an indeterminacy point of $f$ located on a fiber $F_x\subset X$ such that $\overline{f}(x)\neq x$. 
\begin{enumerate}
\item If $\overline{f}$ has the form $x\mapsto x+1$ then $\Cent_b(f)$ is isomorphic to $\ZZ$;
\item if $\overline{f}$ has the form $x\mapsto \alpha x$ then $\Cent_b(f)$ is isomorphic to the product of $\ZZ$ with a finite cyclic group. 
\end{enumerate}
Note that the first case does not occur when $\charK\neq 0$.
\end{cor}
\begin{pro}
Proposition \ref{persjonq} shows that the birational transformation $f$ satisfies the hypothesis of Proposition \ref{orbitargument}. Let $n$ denote the number of persistent indeterminacy points of $f$ with infinite backward orbits. Let $g\in\Cent(f)$. It is in $\Jonq$ by Lemma \ref{jonqcentjonq}. Proposition \ref{orbitargument} says that $g^{n!}\circ f^l$ preserves every member of a pencil of rational curves for some $l\in \ZZ$. The proof of Proposition \ref{orbitargument} shows that infinitely many members of this pencil of rational curves are fibers of the initial rational fibration on $X$. Therefore this pencil of rational curves is the initial rational fibration. This means $\overline{g}^{n!}\circ \overline{f}^l=\Id \in \PGLtK$. 

When $\charK=0$ and $\overline{f}$ is $x\mapsto x+1$, its centralizer in $\PGLtK$ is isomorphic to the additive group $\KK$ and this group is torsion free. Thus, $\Cent_b(f)$ is contained in an infinite cyclic group in which $<\overline{f}>$ has index $\leq n!$. The conclusion follows in this case.

When $\overline{f}$ is $x\mapsto \alpha x$ with $\alpha$ of infinite order, its centralizer in $\PGLtK$ is isomorphic to the multiplicative group $\KK^*$. The difference is that, in this case it is possible that $\overline{g}$ has finite order $\leq n!$. Thus, we may have an additional finite cyclic factor in $\Cent_b(f)$.
\end{pro}

\subsection{Local analysis around a fiber}\label{localanalysis}
Now we need to study the case where there is no persistent indeterminacy points. In this section we will work in the following setting:
\begin{itemize}
  \item Let $f$ be a base-wandering \Jonqui twist. We can suppose that $\overline{f}$ is $x\mapsto \alpha x$ or $x\mapsto x+1$.
	\item Up to taking an algebraically stable model as in Theorem \ref{algstabjonq}, we can suppose that $f$ is a birational transformation of a conic bundle $X$ which satisfies the properties in Theorem \ref{algstabjonq}. We can assume that $X$ is minimal in the sense that $f$ switches the two components of every singular fiber of $X$.
	\item We assume that the only indeterminacy points of $f$ are on the fibers $F_0,F_{\infty}$. This is because these two fibers are the only possible invariant fibers and if $f$ has an indeterminacy point elsewhere then the situation is treated in the previous section.
\end{itemize}
  Without loss of generality, let us suppose that $f$ has an indeterminacy point $p$ on the fiber $F_{\infty}$. By algebraic stability $f^{-1}$ has an indeterminacy point $q\neq p$ on $F_{\infty}$. If $x\in\PP^1$ is not $0$ nor $\infty$, then the orbit of $x$ under $\overline{f}$ is infinite and the fiber $F_x$ is regular. As $f$ has an indeterminacy point on $F_{\infty}$, the fiber $F_{\infty}$ is also regular. Assume that $F_{0}$ is singular, then it is the union of two $(-1)$-curves and $f$ exchanges the two components. Since $\Cent(f)\subset \Cent(f^2)$ and the aim of this section is to prove that $\Cent_b(f)$ is finite by cyclic, it is not harmful to replace $f$ with $f^2$ so that the two components of $F_{0}$ are no more exchanged and we can assume that $F_{0}$ is regular. Thus, we can suppose that \begin{itemize}
		\item the surface $X$ is a Hirzebruch surface.
	\end{itemize}

If $\overline{f}$ is $x\mapsto \alpha x$, then $\Cent_b(f)$ is contained in $\{(x\mapsto \gamma x), \gamma \in\KK^*\}$ and all elements of $\Cent_b(f)$ fix $0$ and $\infty$. Similarly if $\overline{f}$ is $x\mapsto x+1$ then all elements of $\Cent_b(f)$ fix $\infty$. Thus $F_0$ or $F_{\infty}$ is $\Cent(f)$-invariant (under total transforms), we will study the (semi-)local behaviour of the elements in $\Cent(f)$ around such an invariant fiber.

\subsubsection{An infinite chain}\label{infinitechain}
We blow up $X$ at $p,q$ the indeterminacy points of $f,f^{-1}$, obtaining a new surface $X_1$. The fiber of $X_1$ over $0$ is a chain of three rational curves $C_{-1}+C_0+C_1$ where $C_1$ (resp. $C_{-1}$) is the exceptional curve corresponding to $p$ (resp. $q$) and $C_0$ is the strict transform of $F_{\infty}\subset X$. Now $f$ induces a birational transformation $f_1$ of $X_1$. As in Lemma \ref{jonqindlem}, we know that $f_1$ (resp. $f_1^{-1}$) has an indeterminacy point $p_2$ (resp. $q_2$) on $C_1$ (resp. $C_{-1}$) which is disjoint from $C_0$. We then blow up $p_2,q_2$ and repeat the process. We have:
\begin{itemize}
	\item for every $n\in \NN$, a surface $X_n$ on which $f$ induces a birational transformation $f_n$;
	\item the fiber of $X_n$ over $0$ is a chain of rational curves $C_{-n},\cdots,C_0,\cdots,C_n$;
	\item $f_n$ (resp. $f_n^{-1}$) has an indeterminacy point $p_{n+1}$ (resp. $q_{n+1}$) on $C_n$ (resp. $C_{-n}$) disjoint from $C_{n-1}$ (resp. $C_{-(n-1)}$).
\end{itemize}

Let $g$ be a birational transformation of $X$ which commutes with $f$. We already observed that $F_{\infty}$ is an invariant fiber of $g$. If $g$ is regular on $F_{\infty}$, then the commutativity implies that $g$ preserves the set $\{p,q\}$. Suppose that $g$ is not regular on $F_{\infty}$. Then $g$ (resp. $g^{-1}$) has an indeterminacy point $p'$ (resp. $q'$) on $F_{\infty}$. Replacing $g$ by $g^{-1}$ or $f$ by $f^{-1}$, we can suppose that $p'\neq q$. Then for every point $x\in F_{\infty}$ such that $x\neq p, p'$, we have that $g(q)=g(f(x))=f(g(x))$ is a point, thus equals $q$. This further implies $q=q'$. Then we apply the same argument to $g,f^{-1}$, obtaining $p=p'$. In summary, $g$ is either regular on $F_{\infty}$ and preserves $\{p,q\}$, or the set of indeterminacy points of $g, g^{-1}$ on $F_{\infty}$ is exactly $\{p,q\}$. 

We lift $g$ to a birational transformation on $X_n$. By repeating the above arguments, we deduce that for all $n\in\NN$ the two indeterminacy points of $f_n,f_n^{-1}$ on the fiber $F_{\infty}\subset X_n$ coincide with that of $g_n,g_n^{-1}$ if the later exist. This means that for a $C_i$ given, and for sufficiently large $n$, the rational curve $C_i$ is a component of the fiber of $X_n$ and $g_n$ maps it to another component $C_j$ of the fiber. In other words $g$ acts on the infinite chain of rational curves $\sum_{n\in \ZZ}C_n$. The dual graph of this infinite chain of rational curves is a chain of vertices indexed by $\ZZ$. The action of $f$ on the dual graph is just a non trivial translation. The isomorphism group of the dual graph is isomorphic to $\ZZ\rtimes\ZZ/2\ZZ$. Those isomorphisms which commute with a non trivial translation coincide with the subgroup of translations $\ZZ$. The above considerations can be summarized as follows:
\begin{lem}\label{chainaction}
There is a group homomorphism $\Phi:\Cent(f)\rightarrow \ZZ$ such that $g(C_n)=C_{\Phi(g)+n}$ for $g\in\Cent(f)$. An element $g\in\Cent(f)$ is in the kernel of $\Phi$ if and only if $g(C_n)=C_n$ for every $n\in\ZZ$. In other words an element $g$ of the kernel of $\Phi$ is regular on the fiber $F_{\infty}$ and fixes the indeterminacy points of $f,f^{-1}$ on this fiber.
\end{lem}

\begin{lem}\label{nootherind}
Let $g$ be an element of $\Cent(f)$. Let $x\in \PP^1$ be a point not fixed by $\overline{f}$. Then $g$ can not have any indeterminacy points on the fiber $F_x$ over $x$.
\end{lem}
\begin{pro}
By our hypothesis $f$ is regular on all fibers $F_{x_n}$ where $\{x_n,n\in\ZZ\}$ denote the orbit of $x$ under $\overline{f}$. If $g$ had an indeterminacy point $p$ on $F_x$, then $f(p),f^2(p),\cdots$ would give us an infinite number of indeterminacy points of $g$.
\end{pro}

\begin{cor}\label{cyclicprop1}
Suppose that $\overline{f}$ is conjugate to $x\mapsto x+1$ (in particular $\charK=0$). Let $g\in \Cent(f)$ be in the kernel of $\Phi:\Cent(f)\rightarrow \ZZ$. Then $g$ is an automorphism of $X$. Furthermore $g$ preserves the rational fibration fiber by fiber.
\end{cor}
\begin{pro}
Lemma \ref{chainaction} says that $g$ does not have any indeterminacy point on the fiber $F_{\infty}$. Lemma \ref{nootherind} says that $g$ does not have any indeterminacy point elsewhere neither. Thus, $g$ is an automorphism. Since $\overline{g}$ commutes with $\overline{f}:x\mapsto x+1$, $\overline{g}$ is $x\mapsto x+v$ for some $v\in\KK$. Suppose by contradiction that $v\neq 0$. Then $g$ is an elliptic element of infinite order and $f\in\Cent(g)$. We can apply Theorem \ref{ellipticthm} to $g,f$ and put them in normal form. As $f$ is a \Jonqui twist, the rational fibration preserved simultaneously by $f$ and $g$ is unique and it must be the rational fibration appeared in the normal form. Hence, Theorem \ref{ellipticthm} forbids $\overline{f},\overline{g}$ to be both non-trivial and of the form $x\mapsto x+w$ with $w\in \KK$.
\end{pro}

When $\overline{f}$ is of the form $x\mapsto \alpha x$, there are two special fibers $F_0,F_{\infty}$ and the above easy argument does not work.

\subsubsection{Formal considerations along a fiber}
In the rest of this section we will assume that $\overline{f}$ is $x\mapsto \alpha x$. There are two invariant fibers $F_{\infty}$ and $F_0$ in this case. We assume that $f$ has an indeterminacy point $q$ on $F_0$. 

The idea of what we do in the sequel is as follows. Let us look at the case where $\KK=\CC$. The indeterminacy point $q\in F_0$ of $f^{-1}$ is a fixed point of $f$, at which the differential of $f$ has two eigenvalues $0$ and $\alpha$; the fiber directon is superattracting and in the transverse direction $f$ is just $x\mapsto \alpha x$. Therefore there is a local invariant manifold at $q$ for $f$, which is a local holomorphic section of the rational fibration. Likewise, there is a local invariant manifold at $p\in F_0$, the indeterminacy point of $f$. These two local holomorphic sections allow us to conjugate locally holomorphically $f$ to $(\alpha x,a(x) y)$ where $a$ is a germ of holomorphic function. The structure of \Jonqui maps is nice enough to allow us to apply this geometric idea over any field in an elementary way. We need just to work with formal series instead of polynomials. 

From now on we fix $f:(x,y)\dashmapsto (\alpha x, \frac{A(x)y+B(x)}{C(x)y+D(x)})$ where $\alpha\in\KK^*$ has infinite order and $A,B,C,D\in \KK[x]$. Without loss of generality, we suppose that 1) the point $(0,0)$ (resp. $(0,\infty)$) is an indeterminacy point of $f$ (resp. $f^{-1}$); 2) one of the $A,B,C,D$ is not a multiple of $x$. This implies 
\begin{equation}
B(0)=C(0)=D(0)=0, \ A(0)\neq 0. \label{eq:fatzero}
\end{equation}
We will consider $A,B,C,D$ as elements of the ring of formal series $\KK\llbracket x\rrbracket$. We will also view $f$ as an element of the formal \Jonqui group $\operatorname{PGL}_2(\KK\llpar x\rrpar)\rtimes \KK^*$ whose elements are formal expressions of the form $(\mu x, \frac{a(x)y+b(x)}{c(x)y+d(x)})$ where $\mu\in\KK^*$ and $a,b,c,d$ belong to $\KK\llpar x\rrpar$, the fraction field of $\KK\llbracket x\rrbracket$. 

\paragraph{Normal form.}
We want to conjugate $f$ to a formal expression of the form $(\alpha x,\beta(x) y),\beta \in\KK\llpar x\rrpar$ by some formal expression $(x,\frac{E(x)y+F(x)}{G(x)y+H(x)})$ with $E,F,G,H\in \KK\llbracket x\rrbracket$. This amounts to say that we are looking for $E,F,G,H\in \KK\llbracket x\rrbracket$ such that $EF-GH\neq 0$ and
\[
\begin{pmatrix}E(\alpha x)&F(\alpha x)\\G(\alpha x)&H(\alpha x)\end{pmatrix}^{-1} \begin{pmatrix}A(x)&B(x)\\C(x)&D(x)\end{pmatrix}
\begin{pmatrix}E(x)&F(x)\\G(x)&H(x)\end{pmatrix}
\]
is a diagonal matrix. By writing out the explicit expressions of the up-right entry and the down-left entry of this matrix product, we obtain two equations to solve:
\begin{align}
F(x)H(\alpha x)A(x)+H(x)H(\alpha x)B(x)-F(x)F(\alpha x)C(x)-H(x)F(\alpha x)D(x)=0 \label{eq:efghone}\\
-E(x)G(\alpha x)A(x)-G(x)G(\alpha x)B(x)+E(x)E(\alpha x)C(x)+G(x)E(\alpha x)D(x)=0 \label{eq:efghtwo}
\end{align}

We will use minuscules to denote the coefficients of the formal series, e.g.\ $E(x)=\sum_{i\in \NN}e_ix^i$. Let us first look at the constant terms of equations \eqref{eq:efghone}, \eqref{eq:efghtwo}, they give 
\[
-e_0g_0a_0-g_0^2b_0+e_0^2c_0+e_0g_0d_0=0=f_0h_0a_0+h_0^2b_0-f_0f_0c_0-f_0h_0d_0.
\]
Since $b_0=c_0=d_0=0$ and $a_0\neq 0$ (see Equation \eqref{eq:fatzero}), we must have $e_0g_0=f_0h_0=0$. We can choose $f_0=g_0=0$ and $e_0=h_0=1$, this guarantees in particular that our solution will satisfy $EH-FG\neq 0$. 

Remark that the equations \eqref{eq:efghone} and \eqref{eq:efghtwo} involve respectively only $E,G$ and $F,H$, and they have exactly the same form. So it suffices to show the existence of $E,G$ which satisfy equation \eqref{eq:efghone}. The constant term is done, let us look at the $x$ term. This leads to a linear equation in $e_1,g_1$ with coefficients involving $a_0,b_0,c_0,d_0,e_0,g_0$ and $\alpha$. Therefore there exists at least one solution for $e_1,g_1$. Then we turn to the next term and get a linear equation in $e_2,g_2$, and so on. Hence, we can find $E,F,G,H$ which satisfy the desired properties. To sum up, we have:
\begin{lem}\label{EFGHlem}
There exists $E,F,G,H\in \KK\llbracket x\rrbracket$ such that:
\begin{itemize}
	\item $E(0)=H(0)=1$ and $F(0)=G(0)=0$, in particular $\begin{pmatrix}E&F\\G&H\end{pmatrix}\in \operatorname{PGL}_2(\KK\llpar x\rrpar)$;
	\item $(x,\frac{E(x)y+F(x)}{G(x)y+H(x)})$ conjugates $f$ to $(\alpha x, \beta(x)y)$ for some $\beta\in\KK\llpar x\rrpar$;	
\end{itemize}
\end{lem}

\paragraph{Projective line over $\KK\llpar x\rrpar$.}
We call an element of $\PP^1(\KK\llpar x\rrpar)=\KK\llpar x\rrpar\bigcup\{\infty\}$ a formal section. We say a formal section $\theta(x)$ passes through the origin if $\theta(0)=0$.
An element $u=(\mu x, \frac{a(x)y+b(x)}{c(x)y+d(x)})$ of the formal \Jonqui group $\operatorname{PGL}_2(\KK\llpar x\rrpar)\rtimes \KK^*$ acts on $\PP^1(\KK\llpar x\rrpar)$ in the following way:
\begin{align*}
\theta(x) &\mapsto u\cdot\theta(x)=\begin{cases}\infty \quad \text{if} \ c(\mu^{-1}x)\theta(\mu^{-1}x)+d(\mu^{-1}x)=0 \\ \frac{a(\mu^{-1}x)\theta(\mu^{-1}x)+b(\mu^{-1}x)}{c(\mu^{-1}x)\theta(\mu^{-1}x)+d(\mu^{-1}x)}\quad \text{otherwise}\end{cases},\\
\infty &\mapsto \begin{cases}\infty \quad \text{if} \ c=0\\ \frac{a(\mu^{-1}x)}{c(\mu^{-1}x)} \quad \text{if} \ c\neq 0\end{cases}. 
\end{align*}
Geometrically this is saying that a formal section of the rational fibration is sent to another by a formal \Jonqui transformation.
Remark that this action on $\PP^1_{\KK\llpar x\rrpar}$ is not an automorphism of $\KK\llpar x\rrpar$-algebraic variety. In scheme theoretic language, we have a commutative diagram:
\[
\begin{tikzcd}
\PP^1_{\KK\llpar x\rrpar} \arrow{r}{\theta\mapsto u\cdot \theta} \arrow{d}{} & \PP^1_{\KK\llpar x\rrpar}\arrow{d}{}\\
\operatorname{Spec}(\KK\llpar x\rrpar) \arrow{r}{\mu x\mapsfrom x}& \operatorname{Spec}(\KK\llpar x\rrpar).
\end{tikzcd}
\]
Thus, we have a group homomorphism from $\operatorname{PGL}_2(\KK\llpar x\rrpar)\rtimes \KK^*$ to the group of such twisted automorphisms of $\PP^1_{\KK\llpar x\rrpar}$.

Now let $g\in\Cent(f)$ be an element in the kernel of $\Phi$. Recall (see Lemma \ref{chainaction}) that $g$ is regular on the fiber $F_0$ and fixes $(0,0),(0,\infty)$. We showed that $f$ is conjugate by $\begin{pmatrix}E&F\\G&H\end{pmatrix}$ to a formal expression $\hat{f}$ of the form $(\alpha x,\beta(x)y)$. We conjugate $g$ by $\begin{pmatrix}E&F\\G&H\end{pmatrix}$ too to get a formal expression $\hat{g}$. Then $\hat{g}$ commutes with $\hat{f}$. 

Recall that, by Lemma \ref{EFGHlem}, we get $\begin{pmatrix}1&0\\0&1\end{pmatrix}$ when we evaluate the formal expression $\begin{pmatrix}E&F\\G&H\end{pmatrix}$ at $x=0$. Together with the fact that $g\in\Ker(\Phi)$, this implies that we get $y\mapsto \delta_0 y$ for some $\delta_0\in\KK^*$ when we evaluate $\hat{g}$ at $x=0$.

Let us consider the actions of $\hat{f},\hat{g}$ on $\PP^1_{\KK\llpar x\rrpar}$ as described above. Since $\hat{f}$ is in diagonal form, it fixes the points $0$ and $\infty$ of $\PP^1_{\KK\llpar x\rrpar}$. 
\begin{lem}
If $\theta\in\PP^1_{\KK\llpar x\rrpar}$ satisfies $\theta(0)=0$ and $\hat{f}\cdot\theta(x)=\theta(x)$, then $\theta=0$.
\end{lem}
\begin{pro}
The equation $\hat{f}\cdot\theta(x)=\theta(x)$ can be written as $\beta(\alpha^{-1}x)\theta(\alpha^{-1}x)=\theta(x)$, i.e.\ $\theta(\alpha x)^{-1}\beta(x)\theta(x)=1$. Suppose by contradiction that $\theta$ is not $0$. Then we can write $\theta(x)$ as $x^r\tilde{\theta}(x)$ where $r>0$ and $\tilde{\theta}(0)\neq 0$. Hence we have $\tilde{\theta}(\alpha x)^{-1}\beta(x)\tilde{\theta}(x)=\alpha^r$. This implies that $\hat{f}$ is conjugate by $(x,\tilde{\theta}(x)y)$ to $(\alpha x, \alpha^r y)$. Since $\tilde{\theta}(0)\neq 0$ and $\begin{pmatrix}E(0)&F(0)\\G(0)&H(0)\end{pmatrix}=\begin{pmatrix}1&0\\0&1\end{pmatrix}$, this implies that the initial \Jonqui twist $f$ is regular on the fiber $F_0$, contradiction.
\end{pro}

Since $\hat{g}$ is $y\mapsto \delta_0 y$ at $x=0$, it sends the formal section $0\in\PP^1(\KK\llpar x\rrpar)$ to another former section passing through the origin. 
The fact that $\hat{f}$ and $\hat{g}$ commute and the fact that $0$ is the only fixed formal section of $\hat{f}$ which passes through the origin imply that $\hat{g}$ fixes $0\in \PP^1_{\KK\llpar x\rrpar}$. Likewise $\hat{g}$ fixes $\infty$ too. Therefore $\hat{g}$ can be written as $(\gamma x,\delta(x)y)$ where $\gamma\in\KK^*$ and $\delta\in\KK\llpar x\rrpar$ satisfies $\delta(0)=\delta_0\neq 0$.

\paragraph{Normal forms for a pair.}
Let us assume for the moment that $\gamma$ is not a root of unity; we are going to prove that this is impossible. We want to, under this hypothesis, conjugate $\hat{g}=(\gamma x,\delta(x)y)$ to $(\gamma x, \delta(0)y)$ by $h=(x,\xi(x)y)$ for some $\xi\in \KK\llbracket x\rrbracket$. Remark that the conjugate of $\hat{f}$ by $h$ will still be in diagonal form.

We write $\delta=\frac{\omega}{\sigma}$ where $\omega,\sigma\in \KK\llbracket x\rrbracket$ satisfies $\omega(0)\neq 0,\sigma(0)\neq 0$ and $\frac{\omega(0)}{\sigma(0)}=\delta(0)$. We will write $\xi$ as $\sum_{i\in\NN}\xi_ix^i$, and likewise for $\sigma,\omega$.

After conjugation by $h=(x,\xi(x)y)$, $\hat{g}$ becomes 
\[
\tilde{g}=h\circ \hat{g}\circ h^{-1}=(\gamma x,\frac{\xi(\gamma x)}{\xi(x)}\frac{\omega(x)}{\sigma(x)}y).
\]
Therefore the equation we want to solve is
\begin{equation}
\xi(\gamma x)\omega(x)=\frac{\omega_0}{\sigma_0}\xi(x)\sigma(x).
\label{eq:}
\end{equation}
The constant terms of the two sides are automatically equal, let us just choose $\xi_0=1$. Comparing the other terms, we obtain
\begin{align*}
&\xi_0\omega_1+\gamma \xi_1 \omega_0=\frac{\omega_0}{\sigma_0}(\xi_0\sigma_1+\xi_1\sigma_0) \\
&\xi_0\omega_2+\gamma \xi_1\omega_1+\gamma^2\xi_2\omega_0=\frac{\omega_0}{\sigma_0}(\xi_0\sigma_2+\xi_1\sigma_1+\xi_2\sigma_0)\\
&\cdots 
\end{align*}
which are equivalent to 
\begin{align*}
&(\gamma-1)\omega_0\xi_1=\frac{\omega_0}{\sigma_0}\xi_0\sigma_1-\xi_0\omega_1\\
&(\gamma^2-1)\omega_0\xi_2=\frac{\omega_0}{\sigma_0}(\xi_0\sigma_2+\xi_1\sigma_1)-\xi_0\omega_2-\gamma \xi_1\omega_1\\
&\cdots.
\end{align*}
For the $i$-th term, we have a linear equation whose coefficient before $\xi_i$ is $(\gamma^i-1)\omega_0$. Since $\omega\neq 0$ and we have supposed that $\gamma$ is not a root of unity, The above equations always have solutions. In summary, we have the following intermediate lemma (we will get from this lemma a contradiction so its hypothesis is in fact absurd):
\begin{lem}
Suppose that $g\in\Ker(\Phi)$ and the action of $g$ on the base is of infinite order. Then we can conjugate $f$ and $g$, simultaneously by an element in $\operatorname{PGL}_2(\KK\llpar x\rrpar)$ whose evaluation at $x=0$ is $\Id:y\mapsto y$, to 
\[
\tilde{g}=(\gamma x,\delta y),\ \tilde{f}=(\alpha x,\beta(x) y)
\]
where $\alpha,\gamma,\delta\in \KK^*$, $\beta\in \KK\llpar x\rrpar^*$ and $\alpha,\gamma$ are of infinite order.
\end{lem}

Writing down the equation $\tilde{f}\circ\tilde{g}=\tilde{g}\circ\tilde{f}$, we get $\delta\beta(x)=\delta\beta(\gamma x)$. As $\delta\neq 0$, we get $\beta(x)=\beta(\gamma x)$. We write $\beta=\frac{\beta^{num}}{\beta^{den}}$ with $\beta^{num},\beta^{den}\in \KK\llbracket x\rrbracket$ such that at least one of the $\beta^{num}_0,\beta^{den}_0$ is not $0$. The equation becomes
\[
\beta^{num}(x)\beta^{den}(\gamma x)=\beta^{den}(x)\beta^{num}(\gamma x).
\]
By comparing the coeffcients of two sides, we get 
\[
\forall k\in\NN,\ \sum_{i+j=k}{\beta^{num}_i\beta^{den}_j\gamma^j}=\sum_{i+j=k}{\beta^{den}_i\beta^{num}_j\gamma^j}.
\]
Then by induction on $k$ we get from these equations:
\begin{enumerate}
	\item if $\beta^{num}_0=0$ then $\beta^{num}=0$, this is impossible;
	\item if $\beta^{den}_0=0$ then $\beta^{den}=0$), this is again impossible;
	\item if $\beta^{den}\beta^{num}\neq 0$ then $\beta^{num}=\kappa\beta^{den}$ for some $\kappa\in\KK^*$. Then$\tilde{f}=(\alpha x,\kappa y)$, this contradicts the fact that the original birational transformation $f$ has an indeterminacy point on the fiber $F_0$ because to get $\tilde{f}$ we only did conjugations whose evaluation at $x=0$ are the identity $y\mapsto y$.
\end{enumerate}

Thus, we get
\begin{prop}\label{kerphielliptic}
Suppose that $g\in\Ker(\Phi)$. Then $\overline{g}$ has finite order and $g$ is an elliptic element of $\CRK$.
\end{prop}
\begin{pro}
The previous discussion showed that $\overline{g}$ cannot have infinite order. Then an iterate $g^k$ is in $\Jonqz$ and $f\in\Cent(g^k)$. By Theorem \ref{jonqzthm}, an element which commutes with a \Jonqui twist in $\Jonqz$ cannot have an infinite action on the base. As $\overline{f}$ has infinite order, $g^k$ must be elliptic. So $g$ must be elliptic.
\end{pro}

\subsubsection{Another fiber}
The base action $\overline{f}\in\PGLtK$ is $x\mapsto \alpha x$, it has two fixed points $0$ and $\infty$. Recall that we are always under the hypothesis that the indeterminacy points of $f$ are on the fibers $F_0,F_{\infty}$. We have done analysis around the fiber $F_0$ on which $f$ has an indeterminacy point. We will denote by $\Phi_0$ the homomorphism $\Phi$ we considered before. In case $f$ has also an indeterminacy point on $F_{\infty}$, we denote the corresponding homomorphism by $\Phi_{\infty}$. 
\begin{lem}\label{alggpcycliclem}
The image of $\Aut(X)\bigcap\Ker(\Phi_0)\subset\Cent(f)$ in $\Cent_b(f)\subset\PGLtK$ is a finite cyclic group.
\end{lem}
\begin{pro}
We recall first that the automorphism group of a Hirzebruch surface is an algebraic group (see \cite{Mar71}).
An element of $\Cent(f)$ which is regular everywhere on $H$ must be in $\Ker(\Phi_0)$. Thus, $\Aut(X)\bigcap\Ker(\Phi_0)=\Aut(X)\bigcap\Cent(f)$ is an algebraic subgroup of $\Aut(H)$. An automorphism of a Hirzebruch surface always preserves the rational fibration and there is a morphism of algebraic groups from $\Aut(X)$ to $\PGLtK$ (see \cite{Mar71}). The image of $\Aut(X)\bigcap\Ker(\Phi_0)\subset\Cent(f)$ in $\Cent_b(f)\subset\PGLtK$ is an algebraic subgroup $\Lambda$ of $\PGLtK$. By Proposition \ref{kerphielliptic}, the elements of $\Lambda$ are all multiplication by roots of unity. If $\Lambda$ was infinite then it would equal to its Zariski closure in $\PGLtK$ and would be isomorphic to the multiplicative group $\KK^*$. But the existence of a base-wandering \Jonqui twist means that $\KK^*$ contains elements of infinite order, for example $\alpha$. This contradicts the fact that $\Lambda=\KK^*$ is torsion. The conclusion follows.
\end{pro}

We first look at the case where we have two homomorphisms $\Phi_0,\Phi_{\infty}$:
\begin{prop}\label{cyclicprop2}
If $f$ has an indeterminacy point on $F_{\infty}$, then $\Ker(\Phi_0)=\Ker(\Phi_{\infty})$ is a subgroup of $\Aut(X)$. The image of $\Ker(\Phi_0)$ in $\Cent_b(f)\subset\PGLtK$ is a finite cyclic group.
\end{prop}
\begin{pro}
Let $g$ be an element of $\Ker(\Phi_0)$. By Proposition \ref{kerphielliptic} $g$ is an elliptic element of $\CRK$. If $\Phi_{\infty}(g)$ were not trivial, then $g$ would act by a non trivial translation on the corresponding infinite chain of rational curves and could not be conjugate to any automorphism. This means $g$ must belong to $\Ker(\Phi_{\infty})$ and consequently $g$ must be an automorphism of $H$. The second part of the statement follows from Lemma \ref{alggpcycliclem}.
\end{pro}

When $f$ is regular on $F_{\infty}$, we may need to do a little bit more, but we get more precise information as well:
\begin{prop}\label{cyclicprop3}
If $f$ has no indeterminacy points on $F_{\infty}$, then $\Ker(\Phi_0)$ is a finite cyclic group whose elements are automorphisms of $\PP^1\times\PP^1$ of the form $(x,y)\mapsto (\gamma x,y)$ with $\gamma$ a root of unity.
\end{prop}
\begin{pro}
Assume that $f$ is regular on $F_{\infty}$. Let $g\in\Ker(\Phi_0)$ be a non trivial element, it is regular on $F_0$. By Lemma \ref{nootherind}, an indeterminacy point of $g$ can only be located on $F_{\infty}$. 

Suppose that $g$ has an indeterminacy point $p$ on $F_{\infty}$. Then $g^{-1}$ also has an indeterminacy point $q$ on $F_{\infty}$. If $p\neq q$, then $g$ would act by translation on the corresponding infinite chain of rational curves. This means that $g$ would never be conjugate to an automorphism of some surface and contradicts Proposition \ref{kerphielliptic} which asserts that $g$ is elliptic. Thus, we have $p=q$. 
The facts that $f$ commutes with $g$ and that $f$ is regular on $F_{\infty}$ imply $f(p)=p$. We blow up the Hirzebruch surface $X$ at $p$ to get a new surface $X'$ and induced actions $f',g'$. The induced action $f'$ is still regular on the fiber $F_{\infty}'$ and preserves both of the two irreducible components. If $g'$ has an indeterminacy point on $F_{\infty}'$, then as before it coincides with the indeterminacy point of $g'^{-1}$ and must be fixed by $f'$. Then we can keep blowing up indeterminacy points of maps induced from $g$, or contracting $g$-invariant $(-1)$-curves in the fiber, without loosing the regularity of the map induced by $f$. As $g$ is elliptic, we will get at last a surface $\hat{X}$ with induced actions $\hat{f},\hat{g}$ which are both regular on the fiber over $\infty$. We can suppose that $\hat{X}$ is minimal among the surfaces with this property. In particular $\hat{g}$ is an automorphism of $\hat{X}$. Moreover, the proof of Theorem \ref{algstabjonq} shows that $\hat{X}$ is a conic bundle and the only possible singular fiber is $\hat{F}_{\infty}$. We claim that $\hat{F}_{\infty}$ is in fact regular. Suppose by contradiction that $\hat{F}_{\infty}$ is singular. Then it is a chain of two $(-1)$-curves and $\hat{g}$ exchanges the two components. However the conic bundle $\hat{X}$ is obtained from a Hirzebruch surface by a single blow-up, it has a unique section of negative self-intersection which passes through only one of the two components of the singular fiber. As a consequence, the automorphism $\hat{g}$ cannot exchange the two components, contradiction. Thus, replacing $X$ by $\hat{X}$, we can suppose from the beginning that $g$ is an automorphism of the Hirzebruch surface $X$.

Suppose by contradiction that $g$ preserves only finitely many sections of the rational fibration. Since $f$ commutes with $g$, we can assume, after perhaps replacing $f$ by some of its iterates, that $f$ and $g$ preserve simultaneously a section of the rational fibration. Removing this section and the fiber $F_0$ from $H$, we get an open set isomorphic to $\AAA^2$ restricted to which $f$ can be written as $(x',y')\mapsto (\alpha^{-1} x', A(x')y'+B(x'))$ where $A,B\in \KK(x')$. The rational function $A$ must be a constant because $f$ acts as an automorphism on this affine open set. Likewise the rational function $B$ must be a polynomial. But then $(deg(f^n))_{n\in\NN}$ would be a bounded sequence. This contradicts the fact that $f$ is a \Jonqui twist. 

Hence, if $g\in \Ker(\Phi)$ is non-trivial then it preserves necessarily infinitely many sections. This forces $g$ to preserve each member of a pencil of rational curves on $X$ whose general members are sections (see Lemma \ref{invhypersurfaces}). This is only possible if $X=\PP^1\times\PP^1$ and $g$ acts as $(x,y)\mapsto (\gamma x,y)$ with $\gamma\in\KK^*$; here the projection of $\PP^1\times\PP^1$ onto the first factor is the original rational fibration we were looking at. This allows us to conclude by Lemma \ref{alggpcycliclem}.
\end{pro}

\begin{exa}\label{examplecentb}
Let $\mu$ be a $k$-th root of unity, the pair $f:(x,y)\mapsto (\alpha x,\frac{(1+x^k)y+x^k}{(2+x^k)y+1+x^k}), g:(x,y)\mapsto (\mu x, y)$ satisfy the conditions in Proposition \ref{cyclicprop3}.
\end{exa}

Now let $f$ be a base-wandering \Jonqui twist which satisfies the hypothesis made at the beginning of Section \ref{localanalysis}; in particular $f$ is regular outside $F_0\bigcap F_{\infty}$ and $\overline{f}$ is $x\mapsto \alpha x$ or $x\mapsto x+1$.
The image $\Phi_{\infty}(\Cent(f))$ is an infinite cyclic subgroup of $\ZZ$ and is isomorphic to $\ZZ$, it is generated by $\Phi_{\infty}(g)$ for some $g\in\Cent(f)$. Then for any $h\in\Cent(f)$, there exists $k\in\ZZ$ such that $g^{-k}\circ h \in\Ker(\Phi_{\infty})$. Thus, $\overline{g}^{-k}\circ \overline{h}$ belongs to the image of $\Ker(\Phi_{\infty})$ in $\Cent_b(f)$. By Corollary \ref{cyclicprop1}, Proposition \ref{cyclicprop2} and Proposition \ref{cyclicprop3}, the image of $\Ker(\Phi_{\infty})$ in $\Cent_b(f)$ is at worst finite cyclic. Note that $\Cent_b(f)$ is always abelian. Therefore we obtain the last piece of information to prove Theorem \ref{mainthm}: 
\begin{prop}\label{centnonpersjonq}
Let $f$ be a base-wandering \Jonqui twist which satisfies the hypothesis made at the beginning of Section \ref{localanalysis}. Let $g$ be an element of $\Cent(f)$ such that $\Phi_0(g)$ generates the image of $\Phi$. Then $\Cent_b(f)$ is the product of a finite cyclic group with the infinite cyclic group generated by $\overline{g}$.
\end{prop}

\section{Proofs of the main results}\label{proofsection}

\begin{pro}[of Theorem \ref{zzthm}]
The proof is a direct combination of Theorems \ref{loxothm}, \ref{halphenthm}, \ref{ellipticthm}, \ref{jonqzthm} and \ref{mainthm}.
\end{pro}

\begin{pro}[of Theorem \ref{virtuallyabelian}]
Centralizers of loxodromic elements are virtually cyclic by Theorem \ref{loxothm} of Blanc-Cantat. It is proved in \cite{Giz80},\cite{Can11} that centralizers of Halphen twists are virtually abelian (see Theorem \ref{halphenthm}). Centralizers of \Jonqui twists whose actions on the base are of finite order are contained in tori over the function field $\KK(x)$, thus are abelian (\cite{CD12} see Theorem \ref{jonqzthm}). Corollary \ref{basewanderingvirtuallyabelian} says that centralizers of base-wandering \Jonqui twists are virtually abelian when $\charK=0$. Centralizers of infinite order elliptic elements (due to \cite{BD15}) are described in Theorem \ref{ellipticthm}, we use the notations therein. In the first two cases of Theorem \ref{ellipticthm}, by using Remark \ref{pglremark} we see that $\Cent{f}$ is abelian if $\alpha\neq \pm1$, contains an abelian subgroup of index $2$ if $\alpha=-1$ and is not virtually abelian if $\alpha=1$. In the third case of Theorem \ref{ellipticthm}, the kernel of the projection onto the $y$-coordinate is an abelian subgroup of finite index of $\Cent{f}$.
\end{pro}

\begin{pro}[of Remark \ref{degreefunction}]
In the first case the degree function is bounded on $\Gamma$. Indeed there exists $j,k\in\NN^*$ such that $f^j,g^k$ have degree $1$ because $f,g$ are elliptic. Then $f^j,g^k$ generates a subgroup of finite index of $\Gamma$ all of whose non-trivial elements have degree one.

In the second case, the two Halphen twists $f$ and $g$ are automorphisms of a rational surface $X$ preserving an elliptic fibration $X\rightarrow \PP^1$. The elliptic fibration is induced by the linear system corresponding to $mK_X$ for some $m\in\NN^*$. For $n\in\NN$, the actions of $f^n$ and $g^n$ on $Pic(X)$ are respectively
\[
D\mapsto D-mn(D\cdot K_X)\Delta_i+\left( -\frac{m^2}{2}(D\cdot K_X)\cdot (n\Delta_i)^2+m(D\cdot (n\Delta_i)) \right) K_X, \quad i=1,2
\]
where $(\cdot)$ denotes the intersection form and $\Delta_i\in Pic(X)$ satisfies $\Delta_i\cdot K_X=0$ (cf. \cite{Giz80}, \cite{BD15} Section 5). Therefore the action of $f^i\circ g^j$ on $Pic(X)$ is 
\begin{align*}
&D\mapsto D-mi(D\cdot K_X)\Delta_1-mj(D\cdot K_X)\Delta_2+\lambda_{ij} K_X \quad \text{where}\\
& \lambda_{ij}=-\frac{m^2}{2}(D\cdot K_X)\cdot \left (i^2\Delta_1^2 + j^2\Delta_2^2\right) +mD\cdot (i\Delta_1+j\Delta_2)-ijm^2(D \cdot K_X)(\Delta_1\cdot \Delta_2).
\end{align*}
Let $\Lambda$ be an ample class on $X$. Then the degree of $f^i\circ g^j$ is up to a bounded term (cf. \cite{BD15} Section 5)
\begin{align*}
\Lambda\cdot(f^i\circ g^j)^*\Lambda=\Lambda^2-\frac{m^2}{2}(\Lambda\cdot K_X)^2 \left (i^2\Delta_1^2 + j^2\Delta_2^2\right)-ijm^2(\Lambda \cdot K_X)^2(\Delta_1\cdot \Delta_2).
\end{align*}
Note that $\Delta_1^2$ and $\Delta_2^2$ are negative.

The third case is \cite{BC16} Lemma 5.7.


In the fourth case the description of the degree function follows directly from the explicit expressions. 
\end{pro}

\begin{thm}\label{abmaxthm}
Suppose that $\charK=0$. Let $G\subset \CRK$ be a maximal abelian subgroup which has at least one element of infinite order. Then up to conjugation one of the following possibilities holds:
\begin{enumerate}
	\item $G$ is $\{(x,y)\mapsto (\alpha x,\beta y)\vert \alpha,\beta\in \KK^*\}$, $\{(x,y)\mapsto (\alpha x,y+v)\vert \alpha\in \KK^*,v\in \KK\}$ or $\{(x,y)\mapsto (x+u,y+v)\vert u,v\in \KK\}$;
	\item $G$ is the product of $\{(x,y)\mapsto (x,\beta y)\vert \beta\in \KK^*\}$ with an infinite torsion group $G_1$. Each element of $G_1$ has the form 
	\[(x,y)\dashmapsto \left(\eta x,y\frac{S(x)}{S(\eta x)}\right) \quad \eta\ \text{is a root of unity}, S\in\KK(x)\]
	All elements of $G$ are elliptic but $G$ is not conjugate to a group of automorphisms of any rational surface.
	\item $G$ is the product of $\{(x,y)\mapsto (x,y+v)\vert v\in \KK\}$ with an infinite torsion group $G_1$. Each element of $G_1$ has the form 
	\[(x,y)\dashmapsto \left(\eta x,y+S(x)-S(\eta x)\right) \quad \eta\ \text{is a root of unity}, S\in\KK(x)\]
	All elements of $G$ are elliptic but $G$ is not conjugate to a group of automorphisms of any rational surface.
	\item $G$ has a finite index subgroup contained in $\Jonqz=\operatorname{PGL}_2(\KK(x))$.
	\item A finite index subgroup $G'$ of $G$ is a cyclic group generated by a base-wandering \Jonqui twist.
	\item A finite index subgroup $G'$ of $G$ is isomorphic to $\KK^*\times \ZZ$ (resp. $\KK\times \ZZ$) where the first factor is $\{(x,y)\mapsto (x,\beta y)\vert \beta\in\KK^*\}$ (resp. $\{(x,y)\mapsto (x,y+v)\vert v\in\KK\}$) and the second factor is generated by a base-wandering \Jonqui twist, as in the fourth case of Theorem \ref{zzthm};
	\item A finite index subgroup $G'$ of $G$ is isomorphic to $\ZZ^s$ with $s\leq 8$ and $G'$ preserves fiberwise an elliptic fibration;
	\item A finite index subgroup $G'$ of $G$ is a cyclic group generated by a loxodromic element.
\end{enumerate}
\end{thm}
The existences of type two and type three maximal abelian groups are less obvious than the others. We give here two examples.

\begin{exa}\label{torsionexampleone}
Let $q\in\NN^*$. Let $(\xi_n)_n$ be a sequence of elements of $\KK^*$ such that $\xi_n$ is a primitive $q^n$-th root of unity and $\xi_n^q=\xi_{n-1}$ . Let $(R_n)_n$ be a sequence of non-constant rational fractions. For $i\in\NN$, put
\[
f_{i+1}:(x,y)\dashmapsto \left(\xi_{i+1}x,yS_{i+1}(x)\right)  \ \text{with}\ 
S_{i+1}(x)=\frac{R_i(x^{q^i})}{R_i(\xi_1x^{q^i})}\frac{R_{i-1}(x^{q^{i-1}})}{R_{i-1}(\xi_2x^{q^{i-1}})}\cdots\frac{R_1(x)}{R_1(\xi_ix)}.
\]
We have $f_{i+1}^q=f_i$ for all $i\in\NN^*$ so that the group $G_1$ generated by all the $f_i$ is an infinite torsion abelian group. Let $T_i(x)=R_i(x^{q^i})\cdots R_1(x^q)$. The conjugation by $(x,y)\dashmapsto (x,yT_i(x))$ sends the group generated by $f_1,\cdots,f_i$ into the cyclic group $\{(x,y)\mapsto (\xi_i^jx,y)\vert j=0,1,\cdots,q^i-1\}$ whose elements are elliptic. The product of $G_1$ with $\{(x,y)\mapsto (x,\beta y)\vert \beta\in \KK^*\}$ is an abelian subgroup of $\CRK$.

Suppose by contradiction that up to conjugation $G_1$ is a group of automorphisms of a projective rational surface $X$. Since the degree of $f_i$ tends to infinity when $i$ tends to infinity, the $f_i$s are not contained in the identity component $\Aut^0(X)$ for $i$ sufficiently large. Note that the quotient group $\Aut(X)/\Aut^0(X)$ acts on the N\'eron-Severi group $\operatorname{NS}(X,\ZZ)$ with finite kernel (cf. \cite{Har87}, \cite{CanDol12} 1.1) and preserves the intersection form which has signature $(1,\mathrm{dim}\operatorname{NS}(X,\ZZ)-1)$ by Hodge Index Theorem. Therefore infinitely many $f_i$ acts non trivially by isometry on $\operatorname{NS}(X,\ZZ)$. Note that $\Aut^0(X)$ is a finite group because $\Aut(X)/\Aut^0(X)$ is infinite (cf. \cite{Har87}). Since the order of the torsion element $f_i$ tends to infinity, this implies that the order of the action of $f_i$ on $\operatorname{NS}(X,\ZZ)$ tends to infinity. Contradiction because an isometry of $\operatorname{NS}(X,\ZZ)$ cannot have arbitrary high torsion order (the quadratic form is positive definite on the orthogonal of an eigenvector).
\end{exa}

\begin{exa}\label{torsionexampletwo}
We can give an additive version of Example \ref{torsionexampleone}. Let $(\xi_n)_n$ be as in Example \ref{torsionexampleone}. Let $(R_n)_n$ be a sequence of rational fractions whose degrees tend to infinity. For $i\in\NN$, put
\[
f_{i+1}:(x,y)\dashmapsto \left(\xi_{i+1}x,y+S_{i+1}(x)\right)  
\]
with
\[ 
S_{i+1}(x)=R_i(x^{q^i})-R_i(\xi_1x^{q^i})+R_{i-1}(x^{q^{i-1}})-R_{i-1}(\xi_2x^{q^{i-1}})+\cdots+R_1(x)-R_1(\xi_ix).
\]
Let $G_1$ be the group generated by all the $f_i$. The product of $G_1$ with $\{(x,y)\mapsto (x,y+v)\vert v\in \KK\}$ is an abelian subgroup of $\CRK$. Again $G_1$ can not be conjugate to an automorphism group.
\end{exa}


\begin{pro}[of Theorem \ref{abmaxthm}]
Let $G$ be a maximal abelian subgroup of $\CRK$. Note that if $f$ is a non-trivial element of $G$, then $G$ is the maximal abelian subgroup of $\Cent(f)$.

If $G$ contains a loxodromic element $f$, then $G$ is included in $\Cent(f)$ and is virtually the cyclic group generated by $f$ by Theorem \ref{loxothm}; this corresponds to the last case of the above statement. If $G$ contains a Halphen twist, then by Theorem \ref{halphenthm} it is virtually a free abelian group of rank $\leq 8$ which preserves fiberwise an elliptic fibration; this corresponds to the seventh case.

Assume that $G$ contains a base-wandering \Jonqui twist $f$. Theorem \ref{mainthm} says that $\Cent(f)$ is virtually isomorphic to $\KK^*\times \ZZ$, $\KK\times \ZZ$ or $\ZZ$. Thus the same is true for $G$. This correponds to the fifth and the sixth case.

Assume that $G$ contains a non-base-wandering \Jonqui twist $f$. Theorem \ref{jonqzthm} says that $\Cent(f)$ is virtually isomorphic to an abelian subgroup of $\operatorname{PGL}_2(\KK(x))$, so the same is true for $G$. This is the fourth case. 

In the rest of the proof we assume that \emph{$G$ contains only elliptic elements}. Note that $G$ is not necessarily conjugate to a group of automorphisms.

Assume that $G$ contains an element $f:(x,y)\mapsto (\alpha x,y+1)$ with $\alpha\in\KK^*$. By Theorem \ref{ellipticthm} we have
\[
\Cent(f)=\{(x,y)\dashmapsto(\eta(x),y+R(x))\vert \eta\in \PGLtK, \eta(\alpha x)=\alpha \eta(x), R\in \KK(x), R(\alpha x)=R(x).\}
\]
In $\alpha$ has infinite order, then $G=\Cent(f)=\{(x,y)\dashmapsto(\gamma x,y+v)\vert \gamma\in\KK^*,v\in\KK\}$ and we are in the first case. Therefore we can and will assume that $\alpha$ has finite order $d$ so that $f^d$ is $(x,y)\mapsto (x,y+d)$. Assume in this paragraph that $G$ has an element $g$ with an infinite action on the base of the rational fibration $(x,y)\mapsto x$. If the action of $g$ on the base is conjugate to $x\mapsto \beta x$ with $\beta\in\KK^*$, then up to conjugation in $\Jonq$ we can suppose that $g$ is just our initial element $f:(x,y)\dashmapsto (\alpha x,y+1)$ (see Proposition \ref{ellipticnormalform}), so that $G$ is isomorphic to $\KK^*\times \KK$. Now consider the case where the action of $g$ on the base is conjugate to $x'\mapsto x'+1$. The parabolic element $x'\mapsto x'+1$ does not commute with $x\mapsto \alpha x$ (note that $x,x'$ are different coordinates) unless $\alpha=1$. Thus $\alpha=1$ under the existence of such a $g$ and by choosing an appropriate coordinate $x'$, the two elements $g$ and $f$ are respectively $(x',y)\mapsto (x'+1,y+R(x'))$ and $(x',y)\mapsto (x',y+d)$ where $R$ is a polynomial by Lemma \ref{Risapolynomial}. We can conjugate $g$ and $f$, simultaneously by $(x',y)\dashmapsto (x',y+S(x'))$ for some $S\in\KK[X]$, to $(x',y)\mapsto (x'+1,y)$ and $(x',y)\mapsto (x',y+d)$. Then we have
\[
G=\Cent(f^d)\bigcap \Cent(g)=\{(x',y)\mapsto (x'+u,y+v)\vert u,v\in \KK\}.
\] 

We are still under the hypothesis that $G$ contains an element $f:(x,y)\mapsto (\alpha x,y+1)$ with $\alpha\in\KK^*$. Assume now that no element of $G$ has an infinite action on the base of the rational fibration $(x,y)\mapsto x$. Then the description of $\Cent(f)$ implies that $G$ is a subgroup of 
\[\{(x,y)\dashmapsto(\delta x, y+R(x))\vert \delta\in\KK^*, R\in\KK(x)\}.\]
Consider the projection $\pi:G\rightarrow \PGLtK$ which records the action on the base. Denote by $G_0$ the kernel of $\pi$ and by $G_b$ the image of $\pi$. We identify $G_b$ as a subgroup of the multiplicative group of roots of unity of $\KK$. We want to prove that $G_b$ is finite so that $G$ is virtually contained in $\Jonqz=\operatorname{PGL}_2(\KK(x))$. Assume that $G_b$ is an infinite subgroup of the group of roots of unity. We first claim that $G_0$ is isomorphic to $\KK$. Let $h:(x,y)\dashmapsto (x,y+R(x)), R\in\KK(x)$ be an element of $G_0$ and $g:(x,y)\dashmapsto (\beta x,y+S(x)), S\in\KK(x)$ be an element of $G$. The commutation relation $f\circ g=g\circ f$ implies $R(x)=R(\beta x)$. Here $\beta$ can be any element of the infinite group $G_b$. This implies that $R$ is constant, which proves the claim. Let $g:(x,y)\dashmapsto (\gamma x,y+R(x))$ be an element of $G$. Let $d$ be the order of $\gamma$. The transformation $g^d$ is
\[
(x,y)\dashmapsto (x,y+R(x)+R(\gamma x)+\cdots+R(\gamma^{d-1}x))
\]
As $G_0$ is $\KK$, there exists $C\in \KK$ such that $R(x)+\cdots+R(\gamma^{d-1}x)=C$. Define $R^*=R-C/d$. Then we have $R^*(x)+\cdots+R^*(\gamma^{d-1}x)=0$. By additive Hilbert's Theorem 90, there exists $T\in\KK(x)$ such that $R^*(x)=S(\gamma x)-S(x)$. Therefore $R(x)=S(\gamma x)-S(x)+C/d$. Note that $g^*:(x,y)\dashmapsto (\gamma x,y+R^*(x))$ commutes with every element of $G$ once $(x,y)\dashmapsto (\gamma x,y+R^*(x))$ is so. By maximality $g^*$ is also in $G$. Remark that $g^*$ has the same order as $\gamma$. Therefore such $g^*$s form a subgroup of $G$ isomorphic to $G_b$. Thus $G$ is isomorphic to the product $G_0\times G_b$ and we are in the third case.

Assume that $G$ contains an element $f:(x,y)\mapsto(\alpha x,\beta y)$ where $\alpha,\beta\in\KK^*$ and $\beta$ has infinite order.  If $\alpha$ also has infinite order, then Theorem \ref{ellipticthm} implies immediately that $G=\Cent(f)$ is isomorphic to $\KK^*\times \KK^*$ and we are in the first case. Assume that $\alpha$ has finite order but $G$ contains an element $f_1:(x,y)\dashmapsto(\alpha_1 x, yR(x))$ where $R\in\KK(x)$ and $\alpha_1\in\KK^*$ has infinite order. We can suppose that $\alpha=1$ because $\Cent(f)\subset \Cent(f^k)$ for any $k$. We apply Theorem \ref{algstabjonq} to $f_1$ to get a surface $X$ on which the conjugate of $f_1$ acts by automorphism. By Corollary \ref{diagonalisable} $X$ is a Hirzebruch surface. In the proof of Theorem \ref{algstabjonq}, to go from $\PP^1\times \PP^1$ to $X$, we perform elementary transformations at the iterates of the indeterminacy points of $f_1$. Since $f$ commutes with $f_1$ and $f$ preserves every fiber, any iterate of an indeterminacy point of $f_1$ is fixed by $f$. Hence $f$ remains an automorphism on $X$. In other words we can simultaneously conjugate $f,f_1$ to $(x,y)\mapsto(x, \beta y)$ and $(x,y)\mapsto(\alpha_1 x,ry)$ with $r\in\KK^*$. Thus, Theorem \ref{ellipticthm}, when applied respectively to $f$ and $f_1$, shows that $G=\Cent(f)\bigcap \Cent(f_1)$ is isomorphic to the diagonal group $\KK^*\times\KK^*$. Hence we are in the first case. 
 
According to the classification of normal forms of elliptic elements of infinite order (see Proposition \ref{ellipticnormalform}), the only remaining case is the following: $G$ contains an element $f:(x,y)\mapsto(\alpha x,\beta y)$ where $\alpha\in\KK^*$ has finite order and $\beta\in\KK^*$ has infinite order but $G$ contains no elements $(x,y)\dashmapsto(\alpha_1 x, yR(x))$ with $\alpha_1$ of infinite order. In this case $\Cent(f)$ is a subgroup of the \Jonqui group by Theorem \ref{ellipticthm}. Denote by $\pi$ the projection of $G$ into $\PGLtK$. If $\pi(G)$ is finite then we are in the fourth case of Theorem \ref{abmaxthm}. So we assume that $\pi(G)$ is infinite. Then $\pi(G)$ is isomorphic to an infinite subgroup of the group of roots of unity. We want to show that we are in the second case of Theorem \ref{abmaxthm}. By Lemma \ref{telescopic}, each element of $G$ has the form $(x,y)\dashmapsto(\eta(x),y\frac{rS(x)}{S(\eta(x))})$ with $\eta\in\PGLtK,r\in\KK^*,S\in\KK(x)$. If $(x,y)\dashmapsto(\eta(x),y\frac{rS(x)}{S(\eta(x))})$ is an element of $G$ for some $r$, then $(x,y)\dashmapsto(\eta(x),y\frac{S(x)}{S(\eta(x))})$ is also an element of $G$ because it commutes with every other element. However the later has the same order in $G$ as $\eta$ in $\PGLtK$. This means that $G$ has a subgroup isomorphic to $\pi(G)$, so that $G$ is isomorphic to the product of this subgroup with the kernel of $\pi$. To finish the proof, it suffices to show that the kernel of $\pi$ is $\{(x,y)\mapsto (x,\beta y)\vert \beta\in \KK^*\}$. This is because $(x,y)\mapsto (x,\beta y)$ are the only possible elliptic elements by Lemma \ref{telescopic}. 
\end{pro}

\nocite{}

\bibliographystyle{alpha}

\bibliography{biblicentralizers}

\vspace{8mm}
\begin{minipage}{1\textwidth}
 \begin{flushleft}
 ShengYuan Zhao\\
 Institut de Recherche Math\'ematique de Rennes\\
 Universit\'e de Rennes 1\\
 263 avenue du G\'en\'eral Leclerc, CS 74205\\
 F-35042  RENNES C\'edex\\[8mm]
 \vspace{3mm}
 Institute for Mathematical Sciences\\
 Stony Brook University\\
 Stony Brook NY 11794-3660, USA\\
 \emph{e-mail:} \texttt{shengyuan.zhao@stonybrook.edu}\\[8mm]
 \end{flushleft}
\end{minipage}

\end{document}